\documentclass[11pt]{article}
           
\usepackage{amsmath}
\usepackage{amsfonts}
\usepackage{amssymb}
\usepackage[mathscr]{eucal}
              \newif{\ifcomentarios}\comentariosfalse\newtheorem{theorem}{Theorem}

\newtheorem{definition}[theorem]{Definition}

\newtheorem{remark}[theorem]{Remark}

\newenvironment{proof}[1][Proof]{\textit{#1.} }{\hfill $\Box$}

\newcommand{\OBSI}{\begin{remark}\begin{rm}}
\newcommand{\OBSF}{\end{rm}\end{remark}}
\newcommand{\DEFI}{\begin{definition}\begin{rm}}
\newcommand{\DEFF}{\end{rm}\end{definition}}
\newcommand{\zerarcounters}{\setcounter{equation}{0}\setcounter{theorem}{0}}

\newcommand{\dsum}{\displaystyle\sum}

\newcommand{\dint}{\displaystyle\int}

\newcommand{\be}{\begin{eqnarray}}
\newcommand{\en}{\end{eqnarray}}
\newcommand{\bee}{\begin{eqnarray*}}
\newcommand{\ene}{\end{eqnarray*}}

\newcommand{\Z}{\mathbb{Z}}
\newcommand{\N}{\mathbb{N}}
\newcommand{\R}{\mathbb{R}}
\newcommand{\C}{\mathbb{C}}

\newenvironment{proof1}[1][Proof of Theorem~\ref{prop.subor.des}]{\noindent\textit{#1.} }{\hfill $\Box$}
\newenvironment{proof2}[1][Proof of Theorem~\ref{coro.sr}]{\noindent\textit{#1.} }{\hfill $\Box$}

\newenvironment{proof3}[1][Proof of Theorem~\ref{coro.cont.min.cont}]{\noindent\textit{#1.} }{\hfill $\Box$}

\DeclareMathOperator{\dom}{dom}
\DeclareMathOperator{\tr}{Tr}

\DeclareMathOperator{\gerado}{span}
\DeclareMathOperator{\rk}{rank}

\DeclareMathOperator{\hull}{hull}
\DeclareMathOperator{\dist}{dist}

\newtheorem{1}{Definition}[section]

\newtheorem{3}[1]{Lemma}
\newtheorem{4}[1]{Proposition}

\newtheorem{7}[1]{Theorem}

\title{Criteria for the Absolutely Continuous Spectral Components of matrix-valued Jacobi operators}
\author{Fabr\'icio Vieira Oliveira\thanks{fabricio.vieira@engenharia.ufjf.br} \and Silas Luiz de Carvalho\thanks{UFMG}}

\begin{document}
	
\frenchspacing

\maketitle

\begin{abstract}
  We extend in this work the Jitomirskaya-Last inequality~\cite{jito99} and Last-Simon~\cite{simon99} criterion for the absolutely continuous spectral component of a half-line Schr\"odinger operator to the special class of matrix-valued Jacobi operators $H:l^2(\Z,\C^l)\rightarrow l^2(\Z,\C^l)$ given by the law $[H \textbf{u}]_{n} := D_{n - 1} \textbf{u}_{n - 1} + D_{n} \textbf{u}_{n + 1} + V_{n} \textbf{u}_{n}$, where $(D_n)_n$ $(V_n)_n$ are bilateral sequences of $l\times l$ self-adjoint matrices such that $0<\inf_{n\in\Z}s_l[D_n]\le\sup_{n\in\Z}s_1[D_n]<\infty$ (here, $s_k[A]$ stands for the $k$-th singular value of $A$). 

  Moreover, we also show that the absolutely continuous components of even multiplicity of minimal dynamically defined matrix-valued Jacobi operators are constant, extending another result from Last-Simon~\cite{simon99} originally proven for scalar Schr\"odinger operators.
\end{abstract}


\section{Introduction}

A discrete one-dimensional Jacobi operator is defined in $l^{2}(\mathbb{Z}; \mathbb{C})$ by the law
\begin{equation}
\label{eq.schr.dis.esc}
(Hu)_{n} := a_{n}u_{n + 1} + a_{n - 1}u_{n - 1} + v_{n}u_{n},
\end{equation}
where $(a_{n})_{n \in \mathbb{Z}}$ and $(v_{n})_{n \in \mathbb{Z}}$ are bilateral sequences of real numbers, with $a_n>0$ for each $n\in\Z$; the so-called discrete one-dimensional Schr\"odinger operators constitute the particular case where, for each $n\in\Z$, $a_n=1$. Such class of operators has been playing a prominent role in the theory of Mathematical Physics since the advent of Quantum Mechanics. For an account of the main results and applications of the theory of Jacobi operators, see~\cite{teschl00}.

One of the main tools in the study of the spectral properties of such operators consists in relating minimal supports of the spectral types with the asymptotic behavior of the solutions to the eigenvalue equation at $z\in\C$, namely
\begin{equation}
\label{eq.autovalores.schrodinger}
a_{n}u_{n + 1} + a_{n - 1}u_{n - 1} + v_{n}u_{n} = z u_{n}.
\end{equation}

The idea is to associate with such eigenvalue equation, for each $z\in\C$, a bilateral sequence $(\alpha_n(z))_n$ of matrices, the so-called transfer matrices, given by the law
\begin{equation*}
\alpha_{n}(z) := 
\left[
\begin{array}{cc}
\frac{z - v_{n}}{a_n} & -\frac{a_{n-1}}{a_n} \\
1 & 0
\end{array}
\right],
\end{equation*}
and then define another bilateral sequence of matrices 
\begin{equation}
\label{def.mat.trans.erg}
A_{n}(z) := 
\begin{cases}
\alpha_{n - 1}(z)\alpha_{n - 2}(z)...\alpha_{1}(z)\alpha_{0}(\omega), & \mbox{ if }  n \geq 1, \\
& \\
\mathbb{I}_{2}, & \mbox{ if }  n = 0, \\
& \\
\alpha^{-1}_{n}(z)...\alpha^{-1}_{-2}(z)\alpha^{-1}_{-1}(z), & \mbox{ if }   n \leq - 1,
\end{cases}
\end{equation}
the cocycle associated with $(\alpha_n(z))_{n\in\Z}$, so that $(u_{n})_{n \in \mathbb{Z}}$ satisfies \eqref{eq.autovalores.schrodinger} iff 
\begin{equation}
\left[
\begin{array}{c}
u_{n + 1}\\
u_{n}
\end{array}
\right]
=
A_{n}(z)
\left[
\begin{array}{l}
u_{1}\\
u_{0}
\end{array}
\right].
\end{equation}

Then, one may say something about the spectral type at $z\in\C$ by studying the asymptotic properties of $(\Vert A_n(z)\Vert)_n$. In this direction, an important characterization, presented in~\cite{simon99}, relates the absolutely continuous spectral component of the operator \eqref{eq.schr.dis.esc} to the asymptotic behavior of the C\`esaro mean of $(\Vert A_n(z)\Vert)_n$. 
Specifically, it was proved that the set
\begin{equation}
\label{supor.schro}
\mathcal{S} := \{ x \in \mathbb{R}\mid \liminf_{L \rightarrow \infty} \frac{1}{L} \sum^{L}_{n = 1} \left\| A_{n}(x) \right\|^2 < \infty \}
\end{equation}
is a minimal support for the absolutely continuous spectral component of 
$H^+$ (see Subsection~\ref{specsup}), the restriction of the operator $H$ to $l^2(\N;\C)$ (regardless of the boundary condition at $n=0$).

A central result used in the proof of such characterization is the so-called Jitomirskaya-Last inequality, presented in~\cite{jito99}, from which it is possible to obtain many results of the so-called subordinacy theory~\cite{khanpearson,teschl00}. This inequality relates the Weyl-Titchmarsh function $m$, that may be defined as the Borel transform of the spectral measure $\mu$ of $H^+$, namely 
\begin{equation*}
m(z) = \int \frac{1}{x - z} d\mu(x),
\end{equation*}
with the Dirichlet and Neumann solutions, $e$ and $f$, to the eigenvalue equation \eqref{eq.autovalores.schrodinger}, that is, the solutions to \eqref{eq.autovalores.schrodinger} that satisfy, respectively, the initial conditions 
\begin{equation*}
\begin{cases}
e_{0} = & 0, \\
e_{1} = & 1,
\end{cases}
\begin{array}{llll}
& & & \\
& & &  
\end{array}
\begin{cases}
f_{0} = & 1, \\
f_{1} = & 0.
\end{cases}
\end{equation*}
Explicitly, the inequality states that for each $y > 0$, there exists $L(y) \ge 1$ such that 
\begin{equation}
\label{des.jito.schro}
(5 - \sqrt{24}) \frac{\left\|f\right\|_{L(y)}}{\left\|e\right\|_{L(y)}}  < \left| m(x + iy) \right| < (5 + \sqrt{24}) \frac{\left\|f\right\|_{L(y)}}{\left\|e\right\|_{L(y)}},
\end{equation}
where 
\begin{equation*}
\left\|u\right\|_{L} := \left( \sum^{\left\lfloor L \right\rfloor}_{n = 1} \left| u_{n} \right|^{2} + \left(L - \left\lfloor L \right\rfloor\right) \left| u_{\left\lfloor L\right\rfloor + 1} \right|^{2} \right)^{\frac{1}{2}}
\end{equation*}
is the truncated $l^2$-norm of $u$ at $L \ge 1$.

Our main goal in this work is to extend these two results, namely inequality \eqref{des.jito.schro} (a result that is important on its own) and the characterization given by \eqref{supor.schro}, to a particular class of matrix-valued Jacobi operators (see \cite{marx15} for the main results and applications of such operators; see also~\cite{vieiracarvalho21} for an account of the extension of Kotani Theory to ergodic matrix-valued Jacobi operators), defined in $l^{2}(\mathbb{Z}; \mathbb{C}^{l})$ by the law  
\begin{equation} 
\label{eq.ope.din.jacobi}
[H \textbf{u}]_{n} := D_{n - 1} \textbf{u}_{n - 1} + D_{n} \textbf{u}_{n + 1} + V_{n} \textbf{u}_{n},\quad\forall n\in\Z,
\end{equation}
where $(D_{n})_{n \in \mathbb{Z}}$ and $(V_{n})_{n \in \mathbb{Z}}$ are sequences in $M(l, \mathbb{R})$ (the linear space of real $l\times l$ matrices) of self-adjoint matrices such that, for each $n\in\Z$, $D_n$ is invertible. 

Let
\begin{equation}\label{eq.aut.Jacobi}
D_{n - 1} \textbf{u}_{n - 1} + D_{n} \textbf{u}_{n + 1} + V_{n} \textbf{u}_{n}=z\textbf{u}_{n},\quad\forall n\in\Z,
\end{equation}
be the respective eigenvalue equation at $z\in\C$, and define the Dirichlet and Neumann solutions, $\phi$ and $\psi$, as the solutions to~\eqref{eq.aut.Jacobi} that satisfy, respectively, the initial conditions  
\begin{equation}
\label{Neumann.dirich.jacobi}
\begin{cases}
\phi_{0} = & \mathbf{0}, \\
\phi_{1} = & \mathbb{I},
\end{cases}
\begin{array}{llll}
& & & \\
& & &  
\end{array}
\begin{cases}
\psi_{0} = & \mathbb{I}, \\
\psi_{1} = & \mathbf{0}.
\end{cases}
\end{equation}

One may also define the matrix-valued Weyl-Titchmarsh $M$-function of $H^\phi_+:\dom(H^\phi_+)\rightarrow l^2(\Z_+;\C^l)$ (here, $\Z_+:=\N\cup\{0\}$), the restriction of $H$ to $\dom(H^\phi_+):=\{\mathbf{u}\in l^2(\Z_+;\C^l)\mid H\mathbf{u}\in l^2(\Z_+;\C^l),$ $\mathbf{u}_0=0\}$, which satisfies the Dirichlet boundary condition at $n=0$ (such operator is self-adjoint; see Section~\ref{SecReg} for details) as the Borel transform of the matrix-valued spectral measure $(\mu_{\textbf{e}_{1,j}, \textbf{e}_{1,k}}(x))_{j,k}$ associated with $H^\phi_+$, that is, 
\begin{equation}
\label{eq.weyl.matrix}
\left[M^{\phi}(z)\right]_{jk} = \int \frac{1}{x - z} d\mu_{\textbf{e}_{1,j}, \textbf{e}_{1,k}}(x),
\end{equation}
where $\textbf{e}_{1,j}$ is the vector given by 
\begin{equation}
\label{def.vet.can}
(\textbf{e}_{1, j})_{m, n} =
\begin{cases}
1, & \mbox{ if }  m = 1, n = j, \\
& \\
0, & \mbox{otherwise}
\end{cases}
\end{equation}
(see Section~\ref{SecReg} for more details).

More specifically, we prove the following results.

\begin{7}
\label{prop.subor.des}
Let $x \in \mathbb{R}$, $y>0$. Let $M^{\phi}(x + iy)$ be the matrix-valued Weyl-Titchmarsh $M$-function of $H^\phi_+$, given by~\eqref{eq.weyl.matrix}, and let $\phi$, $\psi$ be given by~\eqref{Neumann.dirich.jacobi}. Then, there exist constants $k_{1}, k_{2} > 0$ such that\footnote{The term $\frac{\left\| \phi \right\|^{2}_{L(y)}}{s_{l} \left[\phi\right]^{2}_{L(y)}}$ in the inequality is called condition number. The condition number of a linear map $A$ is the number $\left\| A \right\| \left\| A^{-1} \right\|$.} 
$$
k_{1} \frac{\left\| \psi \right\|_{L(y)}}{\left\| \phi \right\|_{L(y)}} \leq \left\|M^{\phi}(x + iy)\right\| \leq k_{2} \frac{\left\| \psi \right\|_{L(y)}}{\left\| \phi \right\|_{L(y)}}  \frac{\left\| \phi \right\|^{2}_{L(y)}}{s_{l} \left[\phi\right]^{2}_{L(y)}}, 
$$
with $\left\| B \right\|_{L(y)}$ and $s_{l}\left[B\right]_{L}$ denoting, respectively, the truncated matrix norm and the truncated $l$-th singular value of the sequence $B=(B_{n})_{n} \in (M(l, \mathbb{C}))^{\mathbb{N}}$ at $L(y)$ (Definition~\ref{def.norm.trun}), where $L(y)\ge1$ is so that
\begin{equation}
\label{def.l}
2y \left\| D_{0}^{-1} \right\|\left\|\psi \right\|_{L(y)} \left\| \phi \right\|_{L(y)} = 1.
\end{equation}
\end{7}

\begin{remark} The matrix norm considered in Theorem~\ref{prop.subor.des} is the Frobenius norm (see Definition~\ref{frobeniusnorm}).
\end{remark}  

Now, let us assume that the singular values of the $l\times l$ matrices $D_{n}$, $n\in\Z$, satisfy 
\begin{equation}
\label{eq.limit.dn.1}
0 < \inf_{n \in \Z} s_{l}[D_{n}] \leq \sup_{n \in \Z} s_{1}[D_{n}] < \infty; 
\end{equation}
this means, in particular, that $\{\Vert D_n\Vert\}_{n\in\Z}$ and $\{\Vert D_n^{-1}\Vert\}_{n\in\Z}$ are uniformly bounded.

\begin{7}
\label{coro.sr}
Let, for each $r\in\{1,\ldots,l\}$,
\begin{equation}
\label{def.val.r.matriz}
\mathcal{S}_{r} := \{x \in \mathbb{R}\mid \liminf_{L \rightarrow \infty} \frac{1}{L} \sum^{L}_{n = 1} s^{2}_{l - r + 1}[\phi_{n}(x)] + s^{2}_{l - r + 1}[\psi_{n}(x)] < \infty \}.
\end{equation}
Then, the set $\overline{\mathcal{S}_{r} \setminus \mathcal{S}_{r + 1}}^{ess}$ corresponds to the absolutely continuous component of multiplicity $r$ of any self-adjoint extension of the operator $H^+$ restricted to $\dom(H^+_{max}):=\{\mathbf{u}\in l^2(\Z_+;\C^l)\mid H^+\mathbf{u}\in l^2(\Z_+;\C^l)\}$ (satisfying any admissible boundary condition at $n=0$). 
\end{7}


Recall that given a Lebesgue measurable set $B \subseteq \mathbb{R}$, the essential closure of $B$ is the set 
\begin{equation*}
\bar{B}^{ess} := \{x \in \mathbb{R}\mid \kappa((x - \epsilon, x + \epsilon) \cap B) > 0, \forall \epsilon > 0\},
\end{equation*}
where $\kappa$ is the Lebesgue measure.

Besides,  one can rewrite $\mathcal{S}_l$ in~\eqref{def.val.r.matriz} as 
$$
\mathcal{S} := \{x \in \mathbb{R}\mid \liminf_{L \rightarrow \infty}  \frac{1}{L} \sum^{L}_{n = 1}  \left\| A_{n}(x) \right\|^{2} < \infty\},
$$
and so, $\mathcal{S}$ is a minimal support for the absolutely continuous spectral component of multiplicity $l$ of $H^\phi_+$, where $(A_{n}(x))_{n\in\N}$ is the sequence of $2l\times 2l$ transfer matrices associated with the operator \eqref{eq.ope.din.jacobi}, defined analogously to \eqref{def.mat.trans.erg} with, for each $n\in\Z$ and each $z\in\C$,
\begin{equation*}
\alpha_{n}(z) := 
\left[
\begin{array}{cc}
D^{-1}_{n}(z - V_{n}) & -D^{-1}_{n} \\
D_{n} & 0
\end{array}
\right].
\end{equation*}
Thus, $\textbf{u} \in (\mathbb{C}^{l})^{\mathbb{Z}}$ is a solution to the eigenvalue equation~\eqref{eq.aut.Jacobi} at $z \in \mathbb{C}$ iff  
\begin{equation}
\label{equa.auto2}
\left[
\begin{array}{c}
\textbf{u}_{n + 1} \\
D_{n}\textbf{u}_{n}
\end{array}
\right]
=
\alpha_{n}(z)
\left[
\begin{array}{c}
\textbf{u}_{n} \\
D_{n - 1}\textbf{u}_{n - 1}
\end{array}
\right].
\end{equation}

As a consequence of Theorem~\ref{coro.sr}, 
we can also extend another important result presented in \cite{simon99}: namely, we also prove the constancy of the absolutely continuous spectral components of minimal dynamically defined matrix-valued Jacobi operators. 

Let $\Omega \neq \emptyset$ be an arbitrary set, let $T: \Omega \rightarrow \Omega$ be an invertible transformation, let $l \in \mathbb{N}$ and let $D, V: \Omega \rightarrow S(l, \mathbb{R})$ ($S(l, \mathbb{R})$ stands for the linear space of $l\times l$ symmetric and real matrices). A dynamically defined operator is a family of operators in $(\mathbb{C}^{l})^{\mathbb{Z}}$ defined, for each $\omega \in \Omega$, by the law
\begin{equation} 
\label{eq.ope.din}
[H_{\omega} \textbf{u}]_{n} := D(T^{n - 1}\omega) \textbf{u}_{n - 1} + D(T^{n}\omega) \textbf{u}_{n + 1} + V(T^{n}\omega) \textbf{u}_{n},\quad\forall n\in\Z.
\end{equation}

If $\Omega$ is a compact metric space and if $T$ is a minimal transformation (that is, for each $\omega\in\Omega$, $\mathcal{O}(\omega):=\{T^n\omega\mid n\ge 0\}$ is dense in $\Omega$), we prove the following result.

\begin{7}
\label{coro.cont.min.cont}
Let $(H_{\omega})_{\omega \in \Omega}$ be given by \eqref{eq.ope.din} and for each $k\in\{1,\ldots,l\}$, let $\sigma_{ac, 2k}(H_{\omega})$ be the absolutely continuous spectrum of $H_\omega$ of even multiplicity $2k$. Then, for each $\omega_{0}, \omega_{1} \in \Omega$ and each $k\in\{1,\ldots,l\}$, one has 
$$
\sigma_{ac, 2k}(H_{\omega_{0}}) = \sigma_{ac, 2k}(H_{\omega_{1}}),
$$
and so there exist only absolutely continuous components of even multiplicity.
\end{7}


We organize this paper as follows. In Section~\ref{SecReg} we discuss the possible self-adjoint extensions of $H^+$ in the limit point case, and then we define $H^\phi_{\pm}$ and $H^\psi_{\pm}$. We also define the Jost solutions to the eigenvalue equation and present minimal supports for the absolutely continuous spectrum of multiplicity $r\in\{1,\ldots, l\}$ of $H_\pm^\phi$ and $H_\pm^\psi$.

In Section~\ref{subord} we prove Theorem~\ref{prop.subor.des}, in Section~\ref{crit.abs.cont.sp.comp.} we prove Theorem~\ref{coro.sr}, and finally, in Section~\ref{constancy} we prove Theorem~\ref{coro.cont.min.cont}.

In Appendix we present some important results about the Frobenius norm and the singular values of square matrices that are used in Sections~\ref{subord} and~\ref{crit.abs.cont.sp.comp.}.

\section{Regularity of the Operator}
\label{SecReg}

\zerarcounters

\subsection{Self-Adjoint Extensions}

In order to prove our main results, it is required to consider the operators $H^\pm:l^2(\Z_{\pm};\C^l)\rightarrow l^2(\Z_\pm;\C^l)$, 
given by the law
\begin{equation}
\label{def.h.mais}
[H^{\pm} \textbf{u}]_{n} :=
\begin{cases}
\textbf{0}, & \mbox{ if }   n = 0, \\
& \\
D_{n - 1} \textbf{u}_{n - 1} + D_{n} \textbf{u}_{n + 1} + V_{n} \textbf{u}_{n}, & \mbox{ if }  n > 0\; (n<0);
\end{cases}
\end{equation}
here, we set $\Z_+:=\N\cup\{0\}$ and $\Z_-:=\{0,-1,-2,\ldots\}$. Naturally, $H^\pm$ is the restriction of $H$ given by \eqref{eq.ope.din.jacobi} to $l^2(\Z_{\pm};\C^l)$.

We direct our attention to $H^+$, which can be represented as
$$
H^{+} \textbf{u} = 
\left[
\begin{array}{cccccc}
0 & 0 & 0 & 0 &  \ldots \\
D_{0} & V_{1} & D_{1} & 0 & \ldots \\
0 & D_{1} & V_{2} & D_{2} & \ldots \\
0 & 0 & D_{2} & V_{3} & \ldots \\
\vdots & \vdots & \vdots & \vdots & \ddots
\end{array}
\right]
\left[
\begin{array}{c}
\textbf{u}_{0}\\
\textbf{u}_{1}\\
\textbf{u}_{2}\\
\textbf{u}_{3}\\
\vdots
\end{array}
\right].
$$

In general, the operators given by \eqref{eq.ope.din.jacobi} and \eqref{def.h.mais} may be unbounded, and so with domains given by proper subsets of $l^{2}(\mathbb{Z}; \mathbb{C}^{l})$ and $l^{2}(\Z_+; \mathbb{C}^{l})$, respectively. In this case, one needs to determine their respective self-adjoint extensions, if they exist.

Note that $H^{+}$ is not even symmetric; in fact, it is easy to check that, for each $\textbf{x}, \textbf{u} \in (\mathbb{C}^{l})^{\Z_+}$ and each $n \in \Z_+$, one has
\begin{equation}
\label{eq.tripla}
\begin{array}{lll}
\sum^{n}_{k = 1} \left\langle \textbf{x}_{k}, (H^{+}\overline{\textbf{u}})_{k} \right\rangle_{\mathbb{C}^{l}} & = & \sum^{n}_{k = 1} \left\langle (H^{+}\textbf{x})_{k}, \overline{\textbf{u}}_{k} \right\rangle_{\mathbb{C}^{l}} \\
& & \\
& & + \left( \left\langle \textbf{x}_{n}, D_{n}\overline{\textbf{u}}_{n + 1} \right\rangle_{\mathbb{C}^{l}} -  \left\langle \textbf{x}_{n + 1}, D_{n}\overline{\textbf{u}}_{n} \right\rangle_{\mathbb{C}^{l}} \right)\\
& & \\
& & - \left( \left\langle \textbf{x}_{0}, D_{0}\overline{\textbf{u}}_{1} \right\rangle_{\mathbb{C}^{l}} - \left\langle \textbf{x}_{1}, D_{0}\overline{\textbf{u}}_{0} \right\rangle_{\mathbb{C}^{l}} \right);
\end{array}
\end{equation}
therefore,  $H^{+}$ will be a symmetric operator in a domain $\mathcal{D} \subseteq l^{2}(\Z_+, \mathbb{C}^{l})$ iff, for each $\textbf{u}, \textbf{x} \in \mathcal{D}$,
\begin{equation}
\label{eq.wrons.zero}
\lim_{n \rightarrow \infty}
\left( \left\langle \textbf{x}_{n}, D_{n}\overline{\textbf{u}}_{n + 1} \right\rangle_{\mathbb{C}^{l}} -  \left\langle \textbf{x}_{n + 1}, D_{n}\overline{\textbf{u}}_{n} \right\rangle_{\mathbb{C}^{l}} \right) = \left\langle \textbf{x}_{0}, D_{0}\overline{\textbf{u}}_{1} \right\rangle_{\mathbb{C}^{l}} - \left\langle \textbf{x}_{1}, D_{0}\overline{\textbf{u}}_{0} \right\rangle_{\mathbb{C}^{l}}.
\end{equation}

One would like, therefore, to obtain the possible self-adjoint extensions of $H^+$. In order to do that, one may adopt the strategy presented in \cite{carmona90} 
and \cite{teschl00}; one begins with the operators
\begin{eqnarray}
\label{def.oper.h.max}
\begin{array}{llll}H^{+}_{max}: & \dom(H^{+}_{max}) & \rightarrow & l^{2}(\Z_+; \mathbb{C}^{l}),\end{array}
\end{eqnarray}
\begin{eqnarray}
\label{def.oper.h.min}
\begin{array}{llll} H^{+}_{min}: & \dom(H^{+}_{min}) & \rightarrow & l^{2}(\Z_+; \mathbb{C}^{l}),\end{array}
\end{eqnarray}
where 
$$
\dom(H^{+}_{max}):=\{ \textbf{u} \in l^{2}(\Z_+; \mathbb{C}^{l}) \; \vert \; H^{+}\textbf{u} \in l^{2}(\Z_+; \mathbb{C}^{l}) \} 
$$
and 
$$
\dom(H^{+}_{min}):= \{\textbf{x} \in c_{00}(\Z_+;\mathbb{C}^{l})\mid\left\langle \textbf{x}_{0}, D_{0}\overline{\textbf{u}}_{1} \right\rangle_{\mathbb{C}^{l}} - \left\langle \textbf{x}_{1}, D_{0}\overline{\textbf{u}}_{0} \right\rangle_{\mathbb{C}^{l}} = 0, \forall  \textbf{u} \in l^{2}(\Z_+;\mathbb{C}^{l}) \}
$$ 
($c_{00}(\Z_+;\mathbb{C}^{l})$ stands for the linear space of sequences in $\C^l$ with only a finite number of nonzero entries).

\begin{4}
\label{prop.h.min.max}
Let $H^{+}_{max}$ and $H^{+}_{min}$ be defined as above. Then, $H^{+}_{max} = (H^{+}_{min})^{*}$.
\end{4}
\begin{proof}
Firstly, we show that $\dom(H^{+}_{max}) \supseteq \dom((H^{+}_{min})^{*})$. By the definition of adjoint operator, 
$$
\dom((H^{+}_{min})^{*}):=\{\textbf{x} \in l^{2}(\Z_+; \mathbb{C}^{l}) \; \vert \; \exists \textbf{y} \in l^{2}(\Z_+; \mathbb{C}^{l}), \left\langle \textbf{x}, H^{+}_{min}\textbf{u} \right\rangle = \left\langle \textbf{y}, \textbf{u} \right\rangle, \forall \textbf{u} \in \dom(H^{+}_{min})\}. 
$$
Suppose that $\textbf{x} \in \dom((H^{+}_{min})^{*})$. One needs to show 
that $H^{+}\textbf{x} \in l^{2}(\Z_+; \mathbb{C}^{l})$. It follows from relation \eqref{eq.tripla} that for each $\textbf{u} \in \dom(H^{+}_{min})$ and each $n > 1$,
$$
\begin{array}{lll}
 \sum^{n}_{k = 1} \left\langle (H^{+}\textbf{x})_{k}, \overline{\textbf{u}}_{k} \right\rangle_{\mathbb{C}^{l}} & = & \sum^{n}_{k = 1} \left\langle \textbf{x}_{k}, (H^{+}_{min}\overline{\textbf{u}})_{k} \right\rangle_{\mathbb{C}^{l}} \\
& & \\
& & - \left( \left\langle \textbf{x}_{n}, D_{n}\overline{\textbf{u}}_{n + 1} \right\rangle_{\mathbb{C}^{l}} -  \left\langle \textbf{x}_{n + 1}, D_{n}\overline{\textbf{u}}_{n} \right\rangle_{\mathbb{C}^{l}} \right)\\
& & \\
& & + \left(  \left\langle \textbf{x}_{0}, D_{0}\overline{\textbf{u}}_{1} \right\rangle_{\mathbb{C}^{l}} - \left\langle \textbf{x}_{1}, D_{0}\overline{\textbf{u}}_{0} \right\rangle_{\mathbb{C}^{l}} \right).
\end{array}
$$
By letting $n\rightarrow\infty$ on both members of the previous relation and by  the definition of $\dom((H^{+}_{min})^{*})$, one gets 
$$
\sum^{\infty}_{k = 1} \left\langle (H^{+}\textbf{x})_{k}, \overline{\textbf{u}}_{k} \right\rangle_{\mathbb{C}^{l}} = \left\langle \textbf{x}, H^{+}_{min}\textbf{u} \right\rangle = \left\langle \textbf{y}, \textbf{u} \right\rangle.
$$
Then, since $\textbf{u}\in \dom((H^{+}_{min})^{*})$ is arbitrary, it follows that $H^{+}\textbf{x} \in l^{2}(\mathbb{Z}_+; \mathbb{C}^{l})$, and therefore, that $\textbf{x} \in \dom(H^{+}_{max})$.

Suppose now that $\textbf{x} \in \dom(H^{+}_{max})$. One needs to show that for each $\textbf{u} \in \dom(H^{+}_{min})$, $\left\langle \textbf{x}, H^{+}_{min}\textbf{u} \right\rangle = \left\langle H^{+}_{max}\textbf{x}, \textbf{u} \right\rangle$. By relations \eqref{eq.tripla} and \eqref{eq.wrons.zero}, one has 
$$
\sum^{\infty}_{k = 1} \left\langle (H^{+}_{max}\textbf{x})_{k}, \overline{\textbf{u}}_{k} \right\rangle_{\mathbb{C}^{l}} = 
\sum^{\infty}_{k = 1} \left\langle \textbf{x}_{k}, (H^{+}_{min}\overline{\textbf{u}})_{k} \right\rangle_{\mathbb{C}^{l}},
$$
and so, $\textbf{x} \in \dom((H^{+}_{min})^{*})$, proving that $\dom(H^{+}_{max}) \subseteq \dom((H^{+}_{min})^{*})$.
\end{proof}

\begin{remark}\label{Rprop.h.min.max}
  Since $(H^{+}_{min})^{**} = \overline{H^{+}_{min}}$, it follows from Proposition~\ref{prop.h.min.max} that $(H^{+}_{max})^{*} = \overline{H^{+}_{min}}$.
\end{remark}

\

We assume, from now on, that the sequence $(D_{n})_{n}$ satisfies, for each pair $\textbf{u}, \textbf{v} \in l^{2}(\mathbb{Z}_+;\mathbb{C}^{l})$,
\begin{equation}
\label{criterio.ponto}
\lim_{n \rightarrow \infty} \left( \left\langle \textbf{v}_{n}, D_{n}\overline{\textbf{u}}_{n + 1} \right\rangle_{\mathbb{C}^{l}} -  \left\langle \textbf{v}_{n + 1}, D_{n}\overline{\textbf{u}}_{n} \right\rangle_{\mathbb{C}^{l}} \right) = 0
\end{equation}
(this is true, for instance, if $D_{n}$ is uniformly bounded). In this case, the boundary form $\Gamma:(\dom(H^{+}_{min}))^\ast\times(\dom(H^{+}_{min}))^\ast\rightarrow l^2(\Z_+;\C^l)$ associated with $H^{+}_{min}$ is given by 
\begin{equation}
\label{def.gama.front}
\begin{array}{lll}
\Gamma(\textbf{u}, \textbf{v}) & := & \sum^{\infty}_{k = 1} \left\langle ((H^{+}_{min})^\ast\textbf{u})_{k}, \overline{\textbf{v}}_{k} \right\rangle_{\mathbb{C}^{l}} - \left\langle \textbf{u}_{k}, ((H^{+}_{min})^\ast\overline{\textbf{v}})_{k} \right\rangle_{\mathbb{C}^{l}} \\
& & \\
&  = & \left\langle \textbf{u}_{1}, D_{0}\overline{\textbf{v}}_{0} \right\rangle_{\mathbb{C}^{l}} - \left\langle D_0\textbf{u}_{0}, \overline{\textbf{v}}_{1} \right\rangle_{\mathbb{C}^{l}},
\end{array}
\end{equation} 
where we have used, in the last equality, Proposition~\ref{prop.h.min.max}, Remark~\ref{Rprop.h.min.max} and relation~\eqref{eq.wrons.zero}.

Under this hypothesis, one may establish a boundary triple $(l^2(\Z_+;\C^l), \rho_{1}, \rho_{2})$ in order to explicitly obtain the self-adjoint extensions of $H^{+}_{min}$. In what follows, we adopt the strategy presented in~\cite{cesar09}. Define the subspaces 
$$
\begin{array}{lll}
X &:=& \{ \textbf{x} = \textbf{u}_{1} + i D_{0}\textbf{u}_{0}\mid \textbf{u} \in \dom((H^{+}_{min})^{*}) \},  \\
&&  \\
Y &:=& \{ \textbf{y} = \textbf{u}_{1} - i D_{0}\textbf{u}_{0}\mid \textbf{u} \in \dom((H^{+}_{min})^{*}) \},
\end{array}
$$
and the maps $\rho_{1}, \rho_{2}: \dom((H^{+}_{min})^{*}) \rightarrow l^2(\Z_+;\C^l)$, given by
$$
\begin{array}{lll}
\rho_{1}(\textbf{u}) & = & \textbf{u}_{1} + i D_{0}\textbf{u}_{0}, \\
\rho_{2}(\textbf{u}) & = & \textbf{u}_{1} - i D_{0}\textbf{u}_{0}; 
\end{array}
$$
note that $X=\rho_1(\dom(H^{+}_{min})^{*})$ and $Y=\rho_2(\dom(H^{+}_{min})^{*})$.

\begin{3}\label{boundarytriple}
  Let $H^+_{min}$, $\rho_1$ and $\rho_2$ be defined  as above. Then, $(l^2(\Z_+;\C^l), \rho_{1}, \rho_{2})$ is a boundary triple for $H_{min}^+$.
\end{3}
\begin{proof} Let  $\textbf{u}, \textbf{v} \in \dom((H^{+}_{min})^{*})$; it follows from the definitions of $\rho_1$ and $\rho_2$ that
$$
\left\langle \rho_{1}(\textbf{u}), \rho_{1}(\textbf{v}) \right\rangle_{\mathbb{C}^{l}} - \left\langle \rho_{2}(\textbf{u}), \rho_{2}(\textbf{v}) \right\rangle_{\mathbb{C}^{l}} = 2i\Gamma(\textbf{u}, \textbf{v}),
$$
where $\Gamma$ is the boundary form given by \eqref{def.gama.front}.

Moreover, 
we affirm that the deficiency indices of $H^{+}_{min}$, defined as
\[n_{\pm}(H^{+}_{min}):=\dim\ker((H^{+}_{min})^{*}\pm i\mathbb{I}),\]
are equal. In fact, it is known that for symmetric operators, their deficiency indices are constant in each semi-plane of $\C_\pm$ (see page $230$ in~\cite{weid80} for a proof of this statement). If one assumes that the matrices $(D_n)_{n\ge 0}$ and $(V_n)_{n\ge 1}$ are real-valued, the solutions to the eigenvalue equation at $z\in\C$ are conjugated to the solutions to the eigenvalue equation at $\overline{z}$; then, $n_{-}(H^{+}_{min}) = n_{+}(H^{+}_{min})$.

By combining the previous results, one concludes that $(l^2(\Z_+;\C^{l}), \rho_{1}, \rho_{2})$ is a boundary triple for $H_{min}^+$ (see Definition~7.1.11 in~\cite{cesar09}).
\end{proof}

Now, it follows from Lemma~\ref{boundarytriple} and from Theorem~7.1.13 in~\cite{cesar09} that the possible self-adjoint extensions of $H^{+}_{min}$, which will be denoted by $H^+_U$, are precisely
\[\dom(H^+_U):=\{\textbf{u}\in\dom(H^{+}_{min})^\ast\mid\rho_{2}(\textbf{u}) = U \rho_{1}(\textbf{u})\},\qquad  H^+_U(\textbf{u})=(H^{+}_{min})^\ast(\textbf{u}),
\] where $U:Y\rightarrow X$ is an arbitrary unitary map.

Let $\textbf{u}\in\dom(H^+_{min})^\ast$. It follows from the identity $\rho_{2}(\textbf{u}) = U \rho_{1}(\textbf{u})$ that
$$
\textbf{u}_{1} - i D_{0}\textbf{u}_{0} = \rho_{2}(\textbf{u}) = U \rho_{1}(\textbf{u})=U \textbf{u}_{1} + i U D_{0}\textbf{u}_{0},
$$
that is,
$$
(\mathbb{I} - U)\textbf{u}_{1} = i(\mathbb{I} + U)D_{0}\textbf{u}_{0};
$$
thus, if $(\mathbb{I} + U)$ is invertible, one has
$$
i(\mathbb{I} + U)^{-1}(\mathbb{I} - U)\textbf{u}_{1} = -D_{0}\textbf{u}_{0}.
$$

Since $i(\mathbb{I} + U)^{-1}(\mathbb{I} - U)$ is the Cayley transform of $U$, it is a self-adjoint operator (see~\cite{cesar09}). Therefore, in this case, the self-adjoint extensions of $H^+_{min}$ are associated with the self-adjoint operators $i(\mathbb{I} + U)^{-1}(\mathbb{I} - U)$. In particular, if $B:Y\rightarrow Y$ is a self-adjoint operator such that $B \textbf{u}_{1} = -D_{0}\textbf{u}_{0}$ and $B \textbf{v}_{1} = - D_{0}\textbf{v}_{0}$, then
$$
\left\langle \textbf{v}_{0}, D_{0}\textbf{u}_{1} \right\rangle_{\mathbb{C}^{l}} - \left\langle \textbf{v}_{1}, D_{0}\textbf{u}_{0} \right\rangle_{\mathbb{C}^{l}}
=
 \left\langle - B \textbf{v}_{1}, \textbf{u}_{1}\right\rangle_{\mathbb{C}^{l}} - \left\langle \textbf{v}_{1}, - B \textbf{u}_{1} \right\rangle_{\mathbb{C}^{l}}
=
0,
$$ 
and so the boundary form is null.

Under these hypotheses, the self-adjoint extensions of the operator $H^{+}_{min}$ are given by 
\begin{equation}
\label{def.extensoes1}
H_{+}^{B} \textbf{u} = 
\left[
\begin{array}{ccccc}
0 & 0 & 0 & 0 & \ldots \\
0 & (V_{1} - B) & D_{1} & 0 & \ldots \\
0 & D_{1} & V_{2} & D_{2} & \ldots \\
0 & 0 & D_{2} & V_{3} & \ldots \\
\vdots & \vdots & \vdots & \vdots & \ddots
\end{array}
\right]
\left[
\begin{array}{c}
\textbf{u}_{0}\\
\textbf{u}_{1}\\
\textbf{u}_{2}\\
\textbf{u}_{3}\\
\vdots
\end{array}
\right],
\end{equation}
where $B:Y\rightarrow Y$ is any self-adjoint operator that satisfies $B \textbf{u}_{1} = -D_{0}\textbf{u}_{0}$ and $B \textbf{v}_{1} = - D_{0}\textbf{v}_{0}$. In the particular case $B = 0$, the corresponding extension is called Dirichlet operator and it is characterized by the relation $\textbf{u}_{0} = \textbf{0}$. Note that the operator $H$ with domain 
\begin{equation}
\label{eq.dom.h}
\dom H=\{\textbf{x} \in l^{2}(\mathbb{Z}, \mathbb{C}^{l}) \; \vert \; H\textbf{x} \in l^{2}(\mathbb{Z}, \mathbb{C}^{l})\}
\end{equation}
is a finite rank perturbation of the direct sum of the operators $H_{+}^{0}$ and $H_{-}^{0}$, with domains, respectively, in $l^{2}(\mathbb{Z}_{+}, \mathbb{C}^{l})$ and $l^{2}(\mathbb{Z}_{-}, \mathbb{C}^{l})$ ($H_{-}^{0}$ is defined, in $l^{2}(\mathbb{Z}_{-}, \mathbb{C}^{l})$, by a relation analogous to~\ref{def.extensoes1}). 

Naturally, there exist unitary operators $U:Y\rightarrow Y$ such that $(\mathbb{I} + U)$ is singular, and so, 
there exist self-adjoint extensions of $H^+_{min}$ that cannot be represented in the form \eqref{def.extensoes1}. This is the case, for instance, when $U = -\mathbb{I}$, whose corresponding extension is called Neumann operator. More generally, if $(\mathbb{I} - U)$ is invertible, the condition $\rho_{2}(\textbf{u}) = U \rho_{1}(\textbf{u})$ becomes
$$
\textbf{u}_{1} = i(\mathbb{I} - U)^{-1}(\mathbb{I} + U)D_{0}\textbf{u}_{0},
$$
and one can represent such self-adjoint extensions of $H^+_{min}$ as  
\begin{equation}
\label{def.extensoes2}
H_{+}^{C} \textbf{u} = 
\left[
\begin{array}{ccccc}
0 & 0 & 0 & 0 & \ldots \\
(D_{0} - V_{1}CD_{0}) & 0 & D_{1} & 0 & \ldots \\
-D_{1}CD_{0} & 0 & V_{2} & D_{2} & \ldots \\
0 & 0 & D_{2} & V_{3} & \ldots \\
\vdots & \vdots & \vdots & \vdots & \ddots
\end{array}
\right]
\left[
\begin{array}{c}
\textbf{u}_{0}\\
\textbf{u}_{1}\\
\textbf{u}_{2}\\
\textbf{u}_{3}\\
\vdots
\end{array}
\right],
\end{equation}
where $C$ is again a self-adjoint matrix. In particular, the Neumann operator corresponds to $C = 0$ and is characterized by the relation $\textbf{u}_{1} = \textbf{0}$.

\subsection{The limit point case}

It is important to emphasize that 
deficiency indices of $H_{min}^{+}$, $n_-(H_{min}^{+})=n_+(H_{min}^{+})$, may assume values between $0$ and  $2l$. The hypothesis that there exists exactly $l$ solutions to the eigenvalue equation~\eqref{eq.aut.Jacobi} that belong to $l^2(\Z_+;\C^l)$ guarantees that the resolvent operator of $H^+_{min}$ is uniquely represented; one says, in this case, that $H^+_{min}$ is in the limit point case at $+\infty$.

\begin{1}[Limit Point Case]
\label{def.ponto.lim}
The operator $H^{+}_{min}$ is said to be in the limit point case at $+ \infty$ if there exists $z \in \mathbb{C}_{+}$ such that there exist exactly $l$ solutions to the eigenvalue equation~\eqref{eq.aut.Jacobi} that belong to $l^2(\Z_+;\C^l)$. In other words, $H^{+}_{min}$ is in the limit point case at $+ \infty$  if $n_{-}(H^{+}_{min}) = n_{+}(H^{+}_{min}) = l$. 
\end{1}

\begin{remark}
  Actually, given that the deficiency indices of $H^{+}_{min}$ do not depend on $z$ (see~\cite{carmona90} for a proof of this statement), one may conclude that if there exists $z \in \mathbb{C}_{+}$ such that there exist exactly $l$ solutions to the eigenvalue equation~\eqref{eq.aut.Jacobi} that belong to $l^2(\Z_+;\C^l)$, then this is true for each $z\in\C$. Thus, Definition~\ref{def.ponto.lim} does not depend on $z$.
\end{remark}

The notion of limit point case was introduced in the theory of Weyl circles for scalar Sturm-Liouville operators (see~\cite{coddington55} for an account of the theory; in the scalar case, one has $n_{\pm} \in \{0, 1, 2\}$). Namely, if one considers the restriction of the Sturm-Liouville operator defined in $l^2(\N;\C)$ to  $l^2(\{1,2,\ldots,N\};\C)$ and assumes that $n_{\pm} = 1$, it follows that the parametrizations (namely, the Weyl-Titchmarsh $m$-functions) of the $l^2$-solutions to the eigenvalue equation, which constitute a circle in the complex plane, converge to a single point as $N\rightarrow\infty$ (and so, one has the limit point case at $+\infty$). The single $m$-function obtained in the limit is related with the resolvent operator (as we will see below). If $n_{\pm} = 2$, these circles converge to another circle, and operator is in the so-called limit circle case at $+\infty$.    

Within the context of Sturm-Liouville operators, there is in~\cite{atkinson64} an account of the theory of Weyl circles in higher dimensions where such parametrizations are given by unitary linear maps. In \cite{hinton81}, there is a more detailed analysis of such operators, where the case $n_{\pm} < l$ is seen to be related with the rank of the matrix representation of the resolvent operator. For $n_{\pm} > l$, there is more than one matrix representation for the resolvent operator. The situation where $n_{\pm} = 2l$ is analogous to the limit circle case discussed for scalar operators (see~\cite{shi04} for details). 

So, our next goal is to obtain a sufficient condition for $H^+_{min}$ to be in the limit point case at $+\infty$. As we have seen, the condition \eqref{criterio.ponto} is crucial for defining the self-adjoint extensions (analogous conditions were obtained in \cite{chen04} for a class of scalar Jacobi operators; for more general scalar Sturm-Liouville operators, see~\cite{atkinson81,read86,qi04}; for matrix-valued operators, see \cite{irina17,pleijel69}). 

We also emphasize that there is in \cite{alla03} an analysis of the self-adjoint extensions of matrix-valued Jacobi operators of the form \eqref{def.extensoes1} for the case where the deficiency indices are between $0$ and $l$. In this context, the operator is called undetermined if the deficiency indices are equal to $l$ (see \cite{mir98,mir99,mir01} for a discussion involving the characterization of undetermined operators).

It is important to stress out that we are not interested here in the complete characterization of the limit point case of $H^+_{min}$ at $+\infty$, but only in obtaining sufficient conditions for the self-adjoint extensions of $H^+_{min}$ given by \eqref{def.extensoes1} and \eqref{def.extensoes2} to be in the limit point case at $+\infty$. In particular, we would like to guarantee that condition \eqref{criterio.ponto} (which can be seen as a ``boundary condition at $+\infty$'') is satisfied. For scalar Jacobi operators, there is a sufficient condition for the operator to be in the limit point case (see the discussion after Lemma~2.16 in~\cite{teschl00}, for instance). One may extend this result to the matrix-valued case.

\begin{4}
\label{prop.crit.pont.lim}
Let $H^{+}_{min}$ be defined as above. If the sequence $(D_k)_{k\ge 0}$ 
satisfies
\begin{equation}
\label{criterio.soma.limite}
\sum^{\infty}_{k = 0} \frac{1}{\left\| D_{k} \right\|} = \infty,
\end{equation}
then $H^{+}_{min}$ is in the limit point case at $+ \infty$.
\end{4}
\begin{proof}
  Firstly, we prove that $n_{+}(H^{+}_{min}) \geq l$. Let $z\in\C_+$ and consider the system of canonical vectors $\{\textbf{e}_{1, 1},  \textbf{e}_{1, 2}, \ldots  \textbf{e}_{1, l}\}$ given by~\eqref{def.vet.can} 
  (such system is a spectral basis for $H^+_{min}$). Since, for each $n\in\Z$, $D_n$ is invertible, the linear space formed by the solutions to the eigenvalue equation~\eqref{eq.aut.Jacobi} for $H^+_{min}$ at $z\in\C_+$ has dimension $2l$; moreover, since $(H^{+}_{min} - z)^{-1}$ is bounded, it follows that the system 
$$
\{(H^{+}_{min} - z)^{-1}\textbf{e}_{1, 1},  (H^{+}_{min} - z)^{-1}\textbf{e}_{1, 2}, \ldots  (H^{+}_{min} - z)^{-1}\textbf{e}_{1, l}\}
$$
is linearly independent, with each one of its elements belonging to $l^2(\Z_+;\C^l)$; hence, $n_{+}(H^{+}_{min}) \geq l$. 

Now, suppose that $n_{+}(H^{+}_{min}) > l$. 
Then, there exist $z \in \mathbb{C}$ and two linearly independent solutions to the eigenvalue equation, say $\textbf{u}$ and $\textbf{v}$, such that $\textbf{u}, \textbf{v} \in l^{2}(\mathbb{Z}_+; \mathbb{C}^{l})$ and $W_{[\textbf{u}, \textbf{v}]}(1) = C\neq 0$, where $W_{[\textbf{u}, \textbf{v}]}(1)=\left\langle D_{0}\textbf{u}_{1}, \bar{\textbf{v}}_{0} \right\rangle_{\mathbb{C}^{l}} - \left\langle D_{0} \textbf{v}_{1}, \overline{\textbf{u}}_{0}  \right\rangle_{\mathbb{C}^{l}}$ is the Wronskian of $\textbf{u}$ and $\textbf{v}$ at $n=1$. By the constancy of the Wronskian (see Subsection~\ref{resope}), it follows that for each $k\in\Z_+$, 
$$
W_{[\textbf{u}, \textbf{v}]}(k) =  \left\langle D_{k - 1}\textbf{u}_{k}, \bar{\textbf{v}}_{k - 1} \right\rangle_{\mathbb{C}^{l}} - \left\langle D_{k - 1} \textbf{v}_{k}, \overline{\textbf{u}}_{k - 1}  \right\rangle_{\mathbb{C}^{l}} = C,
$$
that is,
$$
\frac{\left| C \right|}{\left\| D_{k - 1} \right\|} \leq  \left( \left\| \textbf{u}_{k} \right\| \left\| \textbf{v}_{k - 1} \right\| +  \left\| \textbf{v}_{k} \right\| \left\| \textbf{u}_{k - 1} \right\| \right).
$$
By summing both members of this inequality from $k=1$ to $k=n$, it follows from Cauchy-Schwarz inequality that
$$
\begin{array}{lll}
\dsum^{n}_{k = 1} \dfrac{\left| C \right|}{\left\| D_{k - 1} \right\|} & \leq & \dsum^{n}_{k = 1} \left\| \textbf{u}_{k} \right\| \left\| \textbf{v}_{k - 1} \right\| + \dsum^{n}_{k = 1} \left\| \textbf{v}_{k} \right\| \left\| \textbf{u}_{k - 1} \right\| \\
& & \\
& \leq & \sqrt{\left(\sum^{n}_{k = 1} \left\| \textbf{u}_{k} \right\|^{2} \right) \left( \sum^{n}_{k = 1} \left\| \textbf{v}_{k - 1} \right\|^{2} \right)} \\
& & \\
& & + \sqrt{\left(\sum^{n}_{k = 1} \left\| \textbf{v}_{k} \right\|^{2} \right) \left( \sum^{n}_{k = 1} \left\| \textbf{u}_{k - 1} \right\|^{2} \right)} \\
& & \\
& \leq 
& 2 \sqrt{\left(\sum^{n + 1}_{k = 0} \left\| \textbf{u}_{k} \right\|^{2} \right) \left( \sum^{n + 1}_{k = 0} \left\| \textbf{v}_{k} \right\|^{2} \right)}. 
\end{array}
$$

Thus, by letting $n \rightarrow \infty$ on both members of the last relation, one gets
$$
|C|\sum^{\infty}_{k = 0} \frac{1}{\left\| D_{k} \right\|} \leq 2\left\|\textbf{u}\right\| \left\|\textbf{v}\right\|   < \infty,
$$
which is absurd. Therefore, $n_{+}(H^{+}_{min}) =l$. 
\end{proof}

In what follows, two of the possible self-adjoint extensions of $H^+_{min}$ will be of particular interest: the Dirichlet operator, given by \eqref{def.extensoes1} and satisfying the boundary condition $\textbf{u}_{0} = \textbf{0}$ in $\dom(H^{+}_{max})$ (the domain of the operator \eqref{def.oper.h.max}), whose matrix form is given by
\begin{equation}
\label{def.oper.diri}
H^{\phi}_{+} \textbf{u} = 
\left[
\begin{array}{ccccc}
0 & 0 & 0 & 0 & \ldots \\
0 & V_{1} & D_{1} & 0 & \ldots \\
0 & D_{1} & V_{2} & D_{2} & \ldots \\
0 & 0 & D_{2} & V_{3} & \ldots \\
\vdots & \vdots & \vdots & \vdots  & \ddots
\end{array}
\right]
\left[
\begin{array}{c}
\textbf{u}_{0}\\
\textbf{u}_{1}\\
\textbf{u}_{2}\\
\textbf{u}_{3}\\
\vdots
\end{array}
\right],
\end{equation} 
and the Neumann operator, given by \eqref{def.extensoes1} and satisfying the boundary condition $\textbf{u}_{1} = \textbf{0}$ in $\dom(H^{+}_{max})$, whose matrix form is given by
\begin{equation}
\label{def.oper.neu}
H^{\psi}_+ \textbf{u} = 
\left[
\begin{array}{ccccc}
0 & 0 & 0 & 0 &  \ldots \\
D_{0} & 0 & D_{1} & 0 &  \ldots \\
0 & 0 & V_{2} & D_{2} & \ldots \\
0 & 0 & D_{2} & V_{3} & \ldots \\
\vdots & \vdots & \vdots & \vdots & \ddots
\end{array}
\right]
\left[
\begin{array}{c}
\textbf{u}_{0}\\
\textbf{u}_{1}\\
\textbf{u}_{2}\\
\textbf{u}_{3}\\
\vdots
\end{array}
\right].
\end{equation}


\subsection{The Resolvent Operator}\label{resope}

Usually, the spectral properties of a self-adjoint operator $T$ are related with the asymptotic behaviour of $(T-z)^{-1}$, with $z=x+iy\in\C_+$, as $y\downarrow 0$. Namely, one may obtain its spectral measures using the Stieltjes inversion formula and the Spectral Theorem. This is particularly true for the Dirichlet and Neumann operators, $H^\phi_+$ and $H_+^\psi$. In order to do that, for each fixed $z\in\rho(H_+^{\phi(\psi)})$, one needs to write the respective resolvent operator in its integral form in terms of the so-called (matrix-valued) Green Function, which is by its turn parametrized by the solutions to the eigenvalue equation~\eqref{eq.aut.Jacobi}.

Let $H^{+}$ be the operator given by \eqref{def.h.mais} and let $z \in \mathbb{C} \setminus \mathbb{R}$. It follows that if $H^+$ is in the limit point case at $+\infty$ (Definition \ref{def.ponto.lim}), then 
$$
\mathcal{J}_{+}(z):=\left\{\textbf{u} \in (\mathbb{C}^{l})^{\mathbb{Z}_{+}} \; \vert \; H^{+}\textbf{u} = z\textbf{u}, \textbf{u}\in l^2(\Z_+;\C^l)\right\}, 
$$
is an $l$-dimension linear subspace.

\begin{1}[Jost Solutions]
\label{def.jost}
Let $H^{+}$ be the operator given by \eqref{def.h.mais} and suppose that it is in the limit point case at $+\infty$. Then, for each $z \in \mathbb{C} \setminus \mathbb{R}$, there exists exactly one sequence $(F^{(+)}_{n}(z))_{n}$ of matrices of size $l \times l$ such that

\begin{eqnarray*}
\begin{array}{ll} (a) & D_{n}F^{(+)}_{n + 1}(z) + D_{n - 1}F^{(+)}_{n - 1}(z) + V_{n} F^{(+)}_{n}(z) = z F^{(+)}_{n}(z), \\
(b) & F^{(+)}_{0}(z) = \mathbb{I}, \\
(c) & \sum^{\infty}_{n = 0} \left\|F^{(+)}_{n}(z) \right\|^{2} < \infty.\end{array}
\end{eqnarray*}
Namely, for each $j\in\{1,\ldots,l\}$, $((F^{(+)}_{n}(z))_j)_{n\ge 0}$ is a solution to the eigenvalue equation~\eqref{eq.aut.Jacobi} for $H^+$ such that $(F^{(+)}_{0}(z))_j=\mathbf{e}_j$, where $(F^{(+)}_{n}(z))_j$ stands for the $j$-th column of $F^{(+)}_{n}(z)$ and $\mathbf{e}_j$ is the $j$-th element of the canonical basis of $\C^l$. Each one of these $l$ solutions is called a Jost solution.
\end{1}

The Jost solutions constitute a basis for $\mathcal{J}_{+}(z)$. For instance, it follows from Proposition~\ref{prop.crit.pont.lim} that if condition~\eqref{criterio.soma.limite} holds, then the operators $H^{\phi}_+$ and $H^{\psi}_+$ are in the limit point case, so there exist exactly $l$ Jost solutions for the respective eigenvalue equations (see also~\cite{kotani88} for a more detailed discussion).

One may also define the sequence of matrices $\left(F^{(-)}_{n}(z)\right)_{n}$ so that
$$
\sum^{ - \infty}_{n = 0} \left\| F^{(-)}_{n}(z) \right\|^{2} < \infty, 
$$
whose sequence of columns are solutions to the eigenvalue equation of $H^-$ at $z$, with the canonical vectors as initial conditions at $n=0$; such solutions are a basis for $\mathcal{J}_{-}(z)$, the space of the solutions to the eigenvalue equation of $H^-$ at $z$ which are square summable at $-\infty$.

In what follows, we present a convenient way to parametrize the spaces $\mathcal{J}_{\pm}(z)$.

\begin{1}[Matrix-valued Weyl-Titchmarsh Functions]
\label{def.m.weyl1}
Let $z \in \mathbb{C} \setminus \mathbb{R}$ and let $\left(F^{(+)}_{n}(z)\right)_{n}$, $\left(F^{(-)}_{n}(z)\right)_{n}$ be as in Definition  \ref{def.jost}. The so-called matrix-valued Weyl-Titchmarsh functions associated with the operators $H_{\pm}^{\phi}$ (given by \eqref{def.oper.diri}) are defined as
\begin{equation}\label{DefMWeyl}
M^{\phi}_{\pm}(z) = - F^{(\pm)}_{\pm 1}(z)D^{-1}_{0}.
\end{equation}
\end{1}

The Weyl-Titchmarsh function $M_{+}^{\phi}(z)$ is also related with the (matrix-valued) spectral measure of the operator $H^{\phi}_{+}$ (see the proof of Proposition~\ref{porp.sup.ac}) by relation \eqref{eq.weyl.matrix}. 

\begin{remark} We note that Definition~\ref{def.m.weyl1} reduces to the definitions of the matrix-valued Weyl-Titchmarsh functions, $M^{\phi}_{\pm}(z)$, presented in~\cite{gesztesy97,kotani88} by taking $D_0=\mathbb{I}$ in \eqref{DefMWeyl}.
  \end{remark}

One may write the Jost solutions as a linear combination of the Dirichlet and Neumann solutions: 
\begin{equation}
\label{equa.m.jost}
F^{(+)}_{n}(z) = \psi_{n}(z) - \phi_{n}(z) M^{\phi}_{+}(z) D_{0};
\end{equation}
this identity is a consequence of the fact that both sides of \eqref{equa.m.jost} are solutions to the eigenvalue equation~\eqref{eq.aut.Jacobi} for $H^\phi_+$ at $z$ that coincide at $n\in\{0,1\}$. 

An important tool in the study of the asymptotic behavior of the solutions to the eigenvalue equation~\eqref{eq.aut.Jacobi} is the so-called Green Formula. Let $\textbf{u}, \textbf{v} \in (\mathbb{C}^{l})^{\mathbb{Z}}$; then, for each integers $n > m$,
\begin{equation}
\label{eq.for.green}
\sum^{n}_{k = m} \left\langle (H\textbf{u})_{k}, \bar{\textbf{v}}_{k} \right\rangle_{\mathbb{C}^{l}} - \left\langle (H\textbf{v})_{k}, \bar{\textbf{u}}_{k} \right\rangle_{\mathbb{C}^{l}} = W_{[\textbf{u}, \textbf{v}]}(n + 1) - W_{[\textbf{u}, \textbf{v}]}(m),
\end{equation}
where $W_{[\textbf{u}, \textbf{v}]}(n)$ is the Wronskian of $\textbf{u}$ and $\textbf{v}$ at $n$, given by
$$
W_{[\textbf{u}, \textbf{v}]}(n) :=  \left\langle D_{n - 1}\textbf{u}_{n}, \bar{\textbf{v}}_{n - 1} \right\rangle_{\mathbb{C}^{l}} - \left\langle D_{n - 1} \textbf{v}_{n}, \overline{\textbf{u}}_{n - 1}  \right\rangle_{\mathbb{C}^{l}}.
$$

If one thinks of $\textbf{u}_{n}$ and $\textbf{v}_{n}$ as column vectors, one gets 
\begin{equation}
\label{wronski.vetor}
W_{[\textbf{u}, \textbf{v}]}(n) =  \textbf{u}^{t}_{n}D_{n - 1}\textbf{v}_{n - 1} - \textbf{v}^{t}_{n}D_{n - 1}\textbf{u}_{n - 1}.
\end{equation}

Moreover, for the case where $(A_{n})$ and $(B_{n})$ are sequences in $M(l,\mathbb{R})$, one obtains, by applying the operator to each of their columns, the following version of Green Formula for matrices:
\begin{equation}
\label{wronski.matriz}
\sum^{n}_{k = m} A^{t}_{k} H(B)_{k}  - H(A)_{k}^{t} B_{k}  =  W_{[A, B]}(n + 1) - W_{[A, B]}(m),
\end{equation}
with
$$
W_{[A, B]}(m) = (A^{t}_{m - 1}D_{m - 1}B_{m} - A^{t}_{m}D_{m - 1}B_{m - 1}).
$$ 

In the specific case that $\textbf{u}_{n}$ and $\textbf{v}_{n}$ are solutions to the same eigenvalue equation, the left side of \eqref{eq.for.green} is zero, from which follows the constancy of the Wronskian; namely, if $\textbf{u}$ and $\textbf{v}$ are solutions to the same eigenvalue equation, then for each $m,n  \in \mathbb{Z}$, 
$$
W_{[\textbf{u}, \textbf{v}]}(n) = W_{[\textbf{u}, \textbf{v}]}(m). 
$$

The next results are extracted from~\cite{vieiracarvalho21}, for which we refer for the proofs.
\begin{4}
\label{prop.m}
Let, for each~$z \in \mathbb{C} \setminus \mathbb{R}$, $\left(F^{(+)}_{n}(z)\right)_{n}$ be as in Definition \ref{def.jost} and let $M_{+}^{\phi}(z)$ be the associated matrix-valued Weyl-Titchmarsh function. Then,
\begin{eqnarray*}
\begin{array}{ll}
(a) & (M_{+}^{\phi}(z))^{t} = M_{+}^{\phi}(z);\\ 
(b) & D_{0}\Im[M_{+}^{\phi}(z)]D_{0} = \Im[z] \sum^{\infty}_{k = 1}  (F^{(+)}_{k}(z))^{*}F^{(+)}_{k}(z).\end{array}
\end{eqnarray*}
\end{4}

\begin{4}[Green Function]
\label{porp.def.green}
Let $H^{\phi}_+$ be given by \eqref{def.oper.diri} and  set, for each $p, q \in \Z_+$ and each $z \in \mathbb{C}\setminus\mathbb{R}$,
\[
G^{\phi}(p, q; z) := 
\left\{
\begin{array}{ll}
- \phi_{p}(z) (D_{0})^{-1} (F^{(+)}_{q}(z))^t, & p \leq q, \\
&\\
- F^{(+)}_{p}(z) (D_{0})^{-1} \phi^{t}_{q}(z), & p > q,
\end{array}\right.
\]
where $F^{(+)}(z)$ and $\phi(z)$ are, respectively, the Jost and Dirichlet solutions to the eigenvalue equation~\eqref{eq.aut.Jacobi} for $H^\phi_+$ at $z$. Then, for each $\textbf{u} \in l^{2}(\Z_+, \mathbb{C}^{l})$,
\begin{equation}
\label{eq.resol.green}
\sum_{q} G^{\phi}(p, q; z)\textbf{u}_{q} = ((H_{+}^{\phi} - z)^{-1}\textbf{u})_{p}.
\end{equation}
\end{4}

\subsection{Spectral Supports}
\label{specsup}

We note that the Green Function $G^{\phi}(1, 1; \cdot):\mathbb{C}_{+}\rightarrow M(l,\mathbb{C})$ is a matrix-valued Herglotz function (that is, $G^{\phi}(1, 1; \cdot)$ is analytic and $\Im G^{\phi}(1, 1;z)>0$, for each $z\in\mathbb{C}_+$, a consequence of the fact that $G$ is the integral kernel of $(H^\phi_+-z)^{-1}$ and $\Im (H^\phi_+-z)^{-1}>0)$, from which follows that for $\kappa$-a.e.~$x\in\mathbb{R}$, 
\[
\lim_{y \downarrow 0}\Im G^{\phi}(1, 1; x \pm iy)<\infty
\] 
(see~\cite{gesztesy97}). By Spectral Theorem one has, for each $\textbf{u} \in l^{2}(\mathbb{Z}_+; \mathbb{C}^{l})$, 
\begin{equation}
\label{eq.med.espec}
\left\langle  (H^{\phi}_+ - z)^{-1} \textbf{u}, \textbf{u} \right\rangle = \int \frac{1}{x - z} d\mu_{\textbf{u}}(x),
\end{equation}
where $\mu_{\textbf{u}}(\cdot):=\langle\textbf{u},E(\cdot)\textbf{u}\rangle$ is a finite Borel measure and $E$ is the resolution of the identity of the operator $H^\phi_+$.

Let $H: \dom(H)\subset\mathcal{H} \rightarrow \mathcal{H}$ be a self-adjoint operator defined in a separable Hilbert space $\mathcal{H}$, and let $\mathcal{B} = \{\textbf{u}_{1}, \textbf{u}_{2}, \ldots, \textbf{u}_{k}\}\subset\mathcal{H}$. The cyclic subspace of $H$ spanned by $\mathcal{B}$ is the space
\begin{equation*}
\mathcal{H}_{\mathcal{B}}:=\overline{\gerado\{\cup_{k=1}^l\{(H)^n(\textbf{u}_{j})\mid n\in\mathbb{N}\}\}}.
\end{equation*}

One says that $\mathcal{B} = \{\textbf{u}_{1}, \textbf{u}_{2}, \ldots, \textbf{u}_{k}\}$ is a spectral basis for $H$ if the system $\mathcal{B}$ is linearly independent and $\mathcal{H}_{\mathcal{B}}=\mathcal{H}$.




In our setting, the $l$ canonical vectors $(\textbf{e}_{1,k})_{k = 1,\ldots,l}$ in $(\mathbb{C}^{l})^{\Z_+}$, where $(\textbf{e}_{1,k})_{n,j} = \delta_{1,n}\delta_{j,k}$, form a spectral basis for $H^{\phi}_+$ (see the proof of Proposition~\ref{prop.crit.pont.lim}). The matrix $(\mu_{\textbf{e}_{1, i}, \textbf{e}_{1, j}})_{1 \leq i,j \leq l}$ is called the spectral (matrix) measure of $H^{\phi}_+$.

The next step consists in obtaining a characterization of the absolutely continuous spectrum of $H^\phi_+$ (including multiplicity) by
establishing minimal supports for the spectral measures.

\begin{1}[Minimal Support]
\label{def.sup.minimal}
One says that a set $S \subseteq \mathbb{R}$ is a minimal support for the positive and finite Borel measure $\mu$ if 
\begin{eqnarray*}\begin{array}{ll}
(i) & \mu (\mathbb{R} \setminus S) = 0; \\
(ii) & S_{0}\subset S,\;\;\; \mu(S_{0}) = 0\qquad \Longrightarrow\qquad \kappa(S_{0}) = 0.
\end{array}\end{eqnarray*} 
\end{1}

In other words, a minimal support for $\mu$ is a Borel set in which $\mu$ is concentrated and such that it is regular with respect to the Lebesgue measure (in the sense that its subsets of zero measure necessarily have zero Lebesgue measure). Definition \ref{def.sup.minimal} induces an equivalence relation in $\mathcal{B}(\R)$ (the Borel $\sigma$-algebra of $\R$): 
$$
S_{1} \sim S_{2} \Longleftrightarrow \kappa(S_{1} \Delta S_{2})  = \mu(S_{1} \Delta S_{2}) = 0,
$$
where $S_{1} \Delta S_{2}:= (S_{1} \setminus S_{2}) \cup (S_{2} \setminus S_{1}) $ is the symmetric difference of $S_{1}$ and $S_{2}$ (see Lemma $2.20$ in \cite{gilbert84} for a proof of this statement).

\begin{4}[Proposition~2.6 in~\cite{vieiracarvalho21}]
\label{porp.sup.ac}
Let $H^{\phi}_+$ be given by \eqref{def.oper.diri}, let $M_{+}^{\phi}$ be the associated matrix-valued Weyl-Titchmarsh function, given by \eqref{def.m.weyl1}, and let $j\in\{1,\ldots,l\}$. Then, 
$$
\Sigma^{\phi}_{ac, j}  :=  \{x \in \mathbb{R}\mid \exists \lim_{y \downarrow 0} M_{+}^{\phi}(x + iy), \rk[\lim_{y \downarrow 0} \Im[M_{+}^{\phi}(x + iy)]] = j\}
$$
is a minimal support for the absolutely continuous component of multiplicity $j$ of the spectral measure, and 
$$
\Sigma^{\phi}_{ac} := \bigcup^{l}_{j = 1} \Sigma^{\phi}_{ac, j}
$$
is a minimal support for the absolutely continuous component. Moreover, 
$$
\Sigma^{\phi}_{s} := \{x \in \mathbb{R}; \lim_{y \downarrow 0} \Im[\tr[M_{+}^{\phi}(x + iy)]] = \infty \}
$$
is a minimal support for the singular component of the spectral measure.
\end{4}
%


Finally, we would like to relate the absolutely continuous spectral components of the operators $H^{\phi,\psi}_+:\dom(H^{+}_{max})\rightarrow l^2(\Z_+;\C^l)$, given by~\eqref{def.oper.diri} and~\eqref{def.oper.neu}, respectively.   

Since 
$$
 H^{\phi}_+ =  H^{\psi}_+ + V_{1} \left( \sum^{l}_{j = 1}  \left\langle . , \textbf{e}_{1, j} \right\rangle  \textbf{e}_{1, j} \right) + D_{1} \left( \sum^{l}_{j = 1}  \left\langle . , \textbf{e}_{1, j} \right\rangle  \textbf{e}_{2, j} \right) - D_{0} \left( \sum^{l}_{j = 1}  \left\langle . , \textbf{e}_{2, j} \right\rangle  \textbf{e}_{0, j}\right),
$$
 %
$H^{\phi}_+$ can be seen as a finite rank perturbation of $H^{\psi}_+$. It is a well-known fact (the so-called Kato-Rosenblum Theorem; see~\cite{simon79}) that a finite rank perturbation of a self-adjoint operator preserves the absolutely continuous components of its spectral measures. 

\begin{4}
\label{prop.supor.pertur}
Let $H_+^{\phi}$ and $H_+^{\psi}$ be the operators given, respectively, by \eqref{def.oper.diri} and \eqref{def.oper.neu}. Then, for each $j\in\{1,\ldots,l\}$, $\sigma_{ac, j}(H^{\phi}_+) = \sigma_{ac, j}(H^{\psi}_+)$.
\end{4}
\begin{proof}
Consider $\{ \textbf{e}_{2, 1}, \textbf{e}_{2, 2}, \ldots, \textbf{e}_{2, l} \}$ as a spectral basis for $H^{\phi}_+$ and $H^{\psi}_+$. Let, for each $j\in\{1,\ldots,l\}$, $H_{j}^{\phi}$ and $H^{\psi}_{j}$ be the restrictions of $H_+^{\phi}$ and $H_+^{\psi}$, respectively, to the subspace spanned by $\textbf{e}_{2, j}$.

Since, for each  $j\in\{1,\ldots,l\}$,
$$
 H^{\phi}_j =  H^{\psi}_j + V_{1}\left\langle . , \textbf{e}_{1, j} \right\rangle  \textbf{e}_{1, j} + D_{1}\left\langle . , \textbf{e}_{1, j} \right\rangle  \textbf{e}_{2, j} - D_{0} \left\langle . , \textbf{e}_{2, j} \right\rangle  \textbf{e}_{0, j},
 $$
 $H_{j}^{\phi}$ is a finite rank perturbation of $H_{j}^{\psi}$; thus, it follows from Kato-Rosenblum Theorem that $H_{j}^{\phi}$ and $H^{\psi}_{j}$ have the same absolutely continuous spectrum.  
\end{proof}

\begin{remark}
Naturally, by using the arguments presented in the proof of Proposition~\ref{prop.supor.pertur}, one can prove that, for each $j\in\{1,\ldots,l\}$, $\sigma_{ac, j}(H^{\phi}_-) = \sigma_{ac, j}(H^{\psi}_-)$.   
\end{remark}

We also want to relate the absolutely continuous components of the spectral measures of the operators $H_{\pm}^{\phi}:\dom(H_{\pm}^{\phi})\rightarrow l^2(\mathbb{Z}_\pm;\C^l)$ given by \eqref{def.oper.diri}, where $\dom(H_\pm^\phi):=\{\textbf{u}\in l^2(\mathbb{Z}_\pm;\C^l)\mid H\textbf{u}\in l^2(\mathbb{Z}_\pm;\C^l),$ $\textbf{u}_0=0\}$, with the absolutely continuous components of the spectral measure of the operator $H:\dom(H)\rightarrow l^2(\mathbb{Z};\C^l)$ given by \eqref{eq.ope.din.jacobi}, where $\dom(H):=\{\textbf{u}\in l^2(\mathbb{Z};\C^l)\mid H\textbf{u}\in l^2(\mathbb{Z};\C^l)\}$. Note that
$$
\begin{array}{lll}
H & = & \left( H_{+}^{\phi} \oplus H_{-}^{\phi} \right) + V_{0}\left( \sum^{l}_{j = 1}  \left\langle . , \textbf{e}_{0, j} \right\rangle  \textbf{e}_{0, j} \right)\\
& & \\
& & +  D_{0}\left( \sum^{l}_{j = 1}  \left\langle . , \textbf{e}_{0, j} \right\rangle  \textbf{e}_{1, j} \right) + D_{0}\left( \sum^{l}_{j = 1}  \left\langle . , \textbf{e}_{1, j} \right\rangle  \textbf{e}_{0, j} \right)\\
& & \\
& & +  D_{-1}\left( \sum^{l}_{j = 1}  \left\langle . , \textbf{e}_{0, j} \right\rangle \textbf{e}_{-1, j} \right) + D_{-1}\left( \sum^{l}_{j = 1}  \left\langle . , \textbf{e}_{-1, j} \right\rangle  \textbf{e}_{0, j} \right). \\ 
\end{array}
$$

The next result relates the absolutely continuous spectral components of $H$ and $H^\phi_-\oplus H^\phi_+$.

\begin{4}[Proposition 2.7 in~\cite{vieiracarvalho21}]
\label{prop.supor.pertur.z}
Let $H$, $H^\phi_-$ and $H^\phi_+$ be as above. 
Then, for each $j\in\{1,\ldots, 2l\}$,  $\sigma_{ac, j}\left( H \right) = \sigma_{ac, j}\left( H_{-}^{\phi} \oplus H_{+}^{\phi} \right)$.
\end{4}

\section{Asymptotic Behavior of Solutions to the Eigenvalue Equation}
\label{subord}

\zerarcounters

One natural way of characterizing the spectral types of a self-adjoint operator $T:\dom(T)\subset\mathcal{H}\rightarrow\mathcal{H}$ is through the study of the  
asymptotic behavior of the solutions to the eigenvalue equation $T\psi=z\psi$,
where $\psi\in\dom(T)$ and $z\in\C$.


In the context of scalar Schr\"odinger operators and, more generally, scalar Sturm-Liouville operators, such analysis can be performed through the so-called subordinacy theory, developed in \cite{khanpearson,gilbert87} (see also~\cite{teschl00}). Furthermore, most of the results obtained in this theory can also be extracted from relation~\eqref{des.jito.schro}, as discussed in \cite{jito99}. 

Theorem~\ref{prop.subor.des} presents a generalization of such relation to matrix-valued Jacobi operators. Note that one needs to consider in this case not only the norm of matrices, but also their singular values. 

Before we present the proof of Theorem~\ref{prop.subor.des}, some preparation is required.

\begin{1}[Truncated Norm]
\label{def.norm.trun}
Let  $L\ge 1$ be a real number, let $k\in\{1,\ldots,l\}$ and let $B=(B_{n})_{n} \in (M(l, \mathbb{C}))^{\mathbb{N}}$. One defines the truncated norm of $B$ at $L$ as 
$$
\left\| 	B \right\|_{L} := \left(    \sum^{\left\lfloor L \right\rfloor}_{n = 1} \left\|B_{n}\right\|^{2} + ( L - \left\lfloor  L \right\rfloor )  \left\|B_{\left\lfloor L\right\rfloor + 1}\right\|^{2} \right)^{\frac{1}{2}},
$$
and its truncated $k$-th singular value at $L$ as
$$
s_{k}\left[ B \right]_{L} := \left( \sum^{\left\lfloor L \right\rfloor}_{n = 1} s^{2}_{k}\left[B_{n}\right] + ( L - \left\lfloor  L \right\rfloor )  s^{2}_{k}\left[B_{\left\lfloor L\right\rfloor + 1}\right] \right)^{\frac{1}{2}}.
$$
\end{1}

One can prove that the truncated norm is, in fact, a norm in $l^2(\N;\C^l)$; in order to prove the triangular inequality, note that
$$
\begin{array}{lll}
\left\| A +	B \right\|_{L} &= &\left( \sum^{\left\lfloor L \right\rfloor}_{n = 1} \left\|A_{n} + B_{n} \right\|^{2} + ( L - \left\lfloor  L \right\rfloor )  \left\|A_{\left\lfloor L\right\rfloor + 1} + B_{\left\lfloor L\right\rfloor + 1}\right\|^{2} \right)^{\frac{1}{2}} \\
&= & \left( \left\Vert A+B\right\Vert^2_{\lfloor L\rfloor} + \left\|( L - \left\lfloor  L \right\rfloor )^{\frac{1}{2}} A_{\left\lfloor L\right\rfloor + 1}+( L - \left\lfloor  L \right\rfloor )^{\frac{1}{2}} B_{\left\lfloor L\right\rfloor + 1}\right\|^{2} \right)^{\frac{1}{2}},
\end{array}
$$
and by considering the triangular inequality for  sequences of size $\left\lfloor  L \right\rfloor$, one gets
\begin{eqnarray*}
\begin{array}{lll}
\left\| A +	B \right\|_{L} & \leq & \left( \left(\Vert A\Vert_{\lfloor L\rfloor}+\Vert B\Vert_{\lfloor L\rfloor}\right)^2 + \left\|( L - \left\lfloor  L \right\rfloor )^{\frac{1}{2}} A_{\left\lfloor L\right\rfloor + 1}+( L - \left\lfloor  L \right\rfloor )^{\frac{1}{2}} B_{\left\lfloor L\right\rfloor + 1}\right\|^{2} \right)^{\frac{1}{2}}\\
&\leq & \left( \Vert A\Vert^2_{\left\lfloor L \right\rfloor}+ \left\|( L - \left\lfloor  L \right\rfloor )^{\frac{1}{2}} A_{\left\lfloor L\right\rfloor + 1}  \right\|^2 \right)^{\frac{1}{2}}+\left( \Vert B\Vert^2_{\left\lfloor L \right\rfloor}+ \left\|( L - \left\lfloor  L \right\rfloor )^{\frac{1}{2}} B_{\left\lfloor L\right\rfloor + 1}  \right\|^2 \right)^{\frac{1}{2}}; 
\end{array}
\end{eqnarray*}
Therefore, 
$$
\left\| A +	B \right\|_{L} \leq \left\| A \right\|_{L} + \left\| B \right\|_{L}.
$$

The main ingredient present in the proof of Theorem~\ref{prop.subor.des} (in analogy to the discussion in \cite{jito99,khanpearson}) consists in relating the truncated norm at $L\ge 1$ of the matrix Jost solutions to the eigenvalue equation at $z= x + iy \in \mathbb{C}_{+}$ with the truncated norm at $L$ of the matrix Neumann and Dirichlet solutions at $x$. 

\begin{3}
\label{lema.1.subor}
Let $z = x + iy \in \mathbb{C}_{+}$, let $\left(F_{n}^{(+)}\right)_{n}$ be as in Definition \ref{def.jost}, and let $\psi_{n}, \phi_{n}$ be, respectively, the Neumann and Dirichlet solutions to the eigenvalue equation~\eqref{eq.aut.Jacobi} for $H^\phi_+$ at $x$. 
Then, for each $n \in \mathbb{Z}_+$, 
$$
\begin{array}{lll} 
F_{n}^{(+)}(z) & = & \psi_{n}(x) - \phi_{n}(x)M^{\phi}_+(z)D_{0} \\
&& \\
&& - iy\psi_{n}(x)\sum^{n}_{k = 1} D_{0}^{-1}\phi^{t}_{k}(x)F_{k}^{(+)}(z) \\
&& \\
&& +  iy\phi_{n}(x)\sum^{n}_{k = 1}D_{0}^{-1}\psi^{t}_{k}(x)F_{k}^{(+)}(z),
\end{array}
$$
where $M^{\phi}_+$ is the matrix-valued Weyl-Titchmarsh function given by Definition \ref{def.m.weyl1}. 
\end{3}

\begin{remark} Although the right-hand side of the previous equation is not defined for $n = 0$, it follows from Definition \ref{def.jost} that $F_{0}^{(+)}(z) = \mathbb{I}$.
\end{remark}

\begin{proof} We omit from the notation the dependence on $x,z$ and $(+)$. Let $(B_n)_{n\ge0}$ be defined, for each $n\in\N$, by the law
$$
B_{n} := \psi_{n} - \phi_{n}M^{\phi}D_{0} - iy\psi_{n}\sum^{n}_{k = 1} D_{0}^{-1}\phi^{t}_{k}F_{k} +  iy\phi_{n}\sum^{n}_{k = 1}D_{0}^{-1}\psi^{t}_{k}F_{k},
$$
with $B_{0} := \mathbb{I}$. In what follows, we show that $B=(B_{n})_{n\ge0}$ satisfies 
the eigenvalue equation \eqref{eq.aut.Jacobi} for $H^\phi_+$ at $z=x+iy$.

Firstly, one may compute separately the terms $D_{n}B_{n + 1}$, $D_{n - 1}B_{n - 1}$ and $(V_{n} - x)B_{n}$. For $D_{n}B_{n + 1}$, one has
$$
\begin{array}{lll}
D_{n}B_{n + 1} & = & D_{n}\left( \psi_{n + 1} - \phi_{n + 1}M^{\phi}D_{0} \right) \\
& & \\
& & - \left( iy D_{n}\psi_{n + 1}\sum^{n + 1}_{k = 1}D_{0}^{-1}\phi^{t}_{k}F_{k} \right) \\
& & \\
& & + \left( iyD_{n}\phi_{n + 1}\sum^{n + 1}_{k = 1}D_{0}^{-1}\psi^{t}_{k}F_{k} \right).
\end{array}
$$
One may write the last two terms as
$$
\begin{array}{lll}
-iy D_{n}\psi_{n + 1}\sum^{n + 1}_{k = 1}D_{0}^{-1}\phi^{t}_{k}F_{k} & = & - iyD_{n}\psi_{n + 1}D_{0}^{-1}\phi^{t}_{n + 1}F_{n + 1} - iyD_{n}\psi_{n + 1}D_{0}^{-1}\phi^{t}_{n}F_{n} \\
& & \\
& & - iy D_{n}\psi_{n + 1}\sum^{n - 1}_{k = 1}D_{0}^{-1}\phi^{t}_{k}F_{k} 
\end{array}
$$
and
$$
\begin{array}{lll}
iyD_{n}\phi_{n + 1}\sum^{n + 1}_{k = 1}D_{0}^{-1}\psi^{t}_{k}F_{k} & = & iyD_{n}\phi_{n + 1}D_{0}^{-1}\psi^{t}_{n + 1}F_{n + 1} + iy D_{n} \phi_{n + 1}D_{0}^{-1}\psi^{t}_{n}F_{n} \\
& & \\
& & + iyD_{n}\phi_{n + 1}\sum^{n - 1}_{k = 1}D_{0}^{-1}\psi^{t}_{k}F_{k}.
\end{array}
$$
For the term $D_{n - 1}B_{n - 1}$, one has
$$
\begin{array}{lll}
D_{n - 1}B_{n - 1} & = & D_{n - 1}\psi_{n - 1} - D_{n - 1}\phi_{n - 1}M^{\phi}D_{0} \\
& & \\ 
& & - iy D_{n - 1}\psi_{n - 1}\sum^{n - 1}_{k = 1}D_{0}^{-1}\phi^{t}_{k}F_{k} \\
& & \\
& &  +  iyD_{n - 1}\phi_{n - 1}\sum^{n - 1}_{k = 1}D_{0}^{-1}\psi^{t}_{k}F_{k}.
\end{array}
$$

At last,
$$
\begin{array}{lll}
(V_{n} - x)B_{n} & = & \left( (V_{n} - x) \psi_{n} - (V_{n} - x) \phi_{n}M^{\phi}D_{0} \right) \\
& & \\
& & - \left( iy (V_{n} - x) \psi_{n}\sum^{n}_{k = 1}D_{0}^{-1}\phi^{t}_{k}F_{k} \right) \\
& & \\
& & + \left( iy (V_{n} - x) \phi_{n}\sum^{n}_{k = 1}D_{0}^{-1}\psi^{t}_{k}F_{k} \right);
\end{array}
$$
again, one can write the last two terms in right-hand side of this identity as
$$
- iy (V_{n} - x) \psi_{n}\sum^{n}_{k = 1}D_{0}^{-1}\phi^{t}_{k}F_{k} = - iy (V_{n} - x) \psi_{n}D_{0}^{-1}\phi^{t}_{n}F_{n}  - iy (V_{n} - x) \psi_{n}\sum^{n - 1}_{k = 1}D_{0}^{-1}\phi^{t}_{k}F_{k}
$$
and
$$
 iy (V_{n} - x) \phi_{n}\sum^{n}_{k = 1}D_{0}^{-1}\psi^{t}_{k}F_{k} = iy (V_{n} - x) \phi_{n}D_{0}^{-1}\psi^{t}_{n}F_{n} + iy (V_{n} - x) \phi_{n}\sum^{n - 1}_{k = 1}D_{0}^{-1}\psi^{t}_{k}F_{k}.
$$

 Finally, when one computes the sum $D_{n}B_{n + 1} + D_{n - 1}B_{n - 1} + (V_{n} - x)B_{n}$, by gathering the terms with the factors $\sum^{n - 1}_{k = 1}D_{0}^{-1}\psi^{t}_{k}F_{k}$ and $\sum^{n - 1}_{k = 1}D_{0}^{-1}\phi^{t}_{k}F_{k}$, one may use the fact 
 that $\phi$ and $\psi$ satisfy the eigenvalue equation~\eqref{eq.aut.Jacobi} for $H^\phi_+$ at $x$ in order to obtain, for each $n\in\N$, 
$$
\begin{array}{lll}
\!\!\!\!\! D_{n}B_{n + 1} + D_{n - 1}B_{n - 1} + (V_{n} - x)B_{n} &= & 
  iyD_{n}\phi_{n + 1}D_{0}^{-1}\psi^{t}_{n + 1}F_{n + 1} - iyD_{n}\psi_{n + 1}D_{0}^{-1}\phi^{t}_{n + 1}F_{n + 1}  \\
&  + & iyD_{n}\phi_{n + 1}D_{0}^{-1}\psi^{t}_{n}F_{n}  -  iyD_{n}\psi_{n + 1}D_{0}^{-1}\phi^{t}_{n}F_{n}  \\
&   + & iy (V_{n} - x) \phi_{n}D_{0}^{-1}\psi^{t}_{n}F_{n}  -  iy (V_{n} - x) \psi_{n}D_{0}^{-1}\phi^{t}_{n}F_{n}.
\end{array}
$$
Now, by rewriting this identity as 
\begin{eqnarray*}
\begin{array}{lll}
D_{n}B_{n + 1} + D_{n - 1}B_{n - 1} + (V_{n} - x)B_{n} &=& iyD_{n} \left( \phi_{n + 1}D_{0}^{-1}\psi^{t}_{n + 1} - \psi_{n + 1}D_{0}^{-1}\phi^{t}_{n + 1}\right)F_{n + 1} \\
  & + & iyD_{n} \left( \phi_{n + 1}D_{0}^{-1}\psi^{t}_{n} - \psi_{n + 1}D_{0}^{-1}\phi^{t}_{n}\right)F_{n}\\
  &+& iy (V_{n} - x) \left( \phi_{n}D_{0}^{-1}\psi^{t}_{n} - \psi_{n}D_{0}^{-1}\phi^{t}_{n}\right)F_{n}.
\end{array}
\end{eqnarray*}
one may conclude that, for each $n\in\N$, 
$$
D_{n}B_{n + 1} + D_{n - 1}B_{n - 1} + (V_{n} - x)B_{n} = iy F_{n}. 
$$

On the other hand, since $(F_{n})$ is a solution to the eigenvalue equation~\eqref{eq.aut.Jacobi} for $H^\phi_+$ at $z=x + iy$, it follows that
$$
D_{n}F_{n + 1} + D_{n - 1}F_{n - 1} + (V_{n} - x)F_{n} = iy F_{n}.
$$

The sequence $(B_{n})$ satisfies, by definition, $B_{0} = \mathbb{I}$ and $B_{1} = - M^{\phi}D_{0}$. One also has from Definition \ref{def.m.weyl1} that $F_{1} = - M^{\phi}D_{0}$,  and from Definition \ref{def.jost} that $F_{0} = \mathbb{I}$. Given that $(F_n)_{n\ge 0}$ and $(B_n)_{n\ge 0}$ are solutions to the eigenvalue equation~\eqref{eq.aut.Jacobi} for $H^\phi_+$ at $z$ such that $F_0=B_0$ and $F_1=B_1$, one concludes that $(F_n)_{n\ge 0}=(B_n)_{n\ge 0}$.
\end{proof}

\

The last ingredient required for the proof of Theorem~\ref{prop.subor.des} consists in showing that for each $y>0$, there exists $L(y)\ge 1$ that satisfies relation~\eqref{def.l}. The $2l$ solutions to the eigenvalue equation~\eqref{eq.aut.Jacobi} associated with the sequences $(\psi_{n})$ and $(\phi_{n})$ span the linear space formed by the solutions to this equation. It follows from the constancy of the Wronskian that the existence of a $l^2$-solution to the eigenvalue  equation at $+\infty$ implies the existence of another solution that is unbounded at $+\infty$. This fact guarantees that $\left\| \psi_{n} \right\|$ or $\left\| \phi_{n} \right\|$ grows indefinitely when $n \rightarrow \infty$. Hence, if one defines the function $f: \mathbb{R}_{+} \rightarrow \mathbb{R}_{+}$ by the law
\begin{equation}
\label{funcaof}  
f(L) := 2\left\| \psi \right\|_{L}\left\| \phi \right\|_L,
\end{equation}
then $f$ is continuous, nondecreasing and such that $\lim_{L \rightarrow \infty}f(L) = \infty$. Therefore, for each $y > 0$ there exists $L=L(y)\ge1$ which satisfies relation~\eqref{def.l}.


\

\begin{proof1}
We omit from the notation the dependence on $z$, $x$ and $(+)$. We also use, throughout the proof, the fact that the function $f$ defined by~\eqref{funcaof} is nondecreasing. It follows from Lemma \ref{lema.1.subor} that for each $n \in \mathbb{N}$,
\begin{equation}
\label{eq.f.alpha}
F_{n} = \psi_{n} - \phi_{n}M^{\phi}D_{0} - iy\psi_{n}\alpha_{n} +  iy\phi_{n}\beta_{n},
\end{equation}
with
\begin{eqnarray*}
\left\{
\begin{array}{lll}
\alpha_{n} & := & \sum^{n}_{k = 1} D_{0}^{-1}\phi^{t}_{k}F_{k},\\
\beta_{n} & := & \sum^{n}_{k = 1}D_{0}^{-1}\psi^{t}_{k}F_{k}.
\end{array}\right.
\end{eqnarray*}

The main idea is to relate the sequences $\psi=(\psi_{n})$ and $\phi=(\phi_{n})_n$ with the matrix $M^{\phi}$ by using Proposition  \ref{prop.m}$-b$ and relation~\eqref{eq.f.alpha}. The first step is to estimate, for each $L\ge 1$, the truncated norm of $\psi\alpha$ and $\phi\beta$ (with the Frobenius matrix norm; see Definition~\ref{frobeniusnorm}), and then use these estimates in order to evaluate $\left\| F \right\|_{L}$. One has
\begin{eqnarray*}
\begin{array}{lll}
\left\| \psi\alpha \right\|^{2}_{L} &= &
\sum^{\left\lfloor  L \right\rfloor}_{n = 1} \tr \left[\left(\sum^{n}_{k = 1} \psi_{n}D_{0}^{-1}\phi^{*}_{k}F_{k}\right)^{*}\left(\sum^{n}_{k = 1} \psi_{n}D_{0}^{-1}\phi^{*}_{k}F_{k}\right)\right] \\
&&\\
& + & ( L - \left\lfloor  L \right\rfloor )  \tr\left[ \left(\sum^{\left\lfloor  L \right\rfloor + 1}_{k = 1}\psi_{\left\lfloor  L \right\rfloor + 1}D_{0}^{-1}\phi^{*}_{k}F_{k}\right)^{*}\left(\sum^{\left\lfloor  L \right\rfloor + 1}_{k = 1}\psi_{\left\lfloor  L \right\rfloor + 1}D_{0}^{-1}\phi^{*}_{k}F_{k}\right)\right] \\
&&\\
&= & \sum^{\left\lfloor  L \right\rfloor}_{n = 1} \tr \left[\left(\sum^{n}_{k = 1} F^{*}_{k}\phi_{k}\right)D_{0}^{-1}\psi^{*}_{n}\psi_{n} D_{0}^{-1}\left(\sum^{n}_{k = 1} \phi^{*}_{k}F_{k}\right)\right] \\
&&\\
&  +& ( L - \left\lfloor  L \right\rfloor )  \tr\left[ \left(\sum^{\left\lfloor  L \right\rfloor + 1}_{k = 1}F^{*}_{k}\phi_{k}\right) D_{0}^{-1}\psi^{*}_{\left\lfloor  L \right\rfloor + 1}\psi_{\left\lfloor  L \right\rfloor + 1}D_{0}^{-1}\left(\sum^{\left\lfloor  L \right\rfloor + 1}_{k = 1}\phi^{*}_{k}F_{k}\right)\right]. 
\end{array}
\end{eqnarray*}

It follows from the commutative property of the trace function that for each $n \in \mathbb{N}$,
$$
\begin{array}{lcl}
& \tr \left[\left(\sum^{n}_{k = 1} F^{*}_{k}\phi_{k}\right)D_{0}^{-1}\psi^{*}_{n}\psi_{n} D_{0}^{-1}\left(\sum^{n}_{k = 1} \phi^{*}_{k}F_{k}\right)\right] &  \\ 
& & \\
=& \tr \left[\psi_{n} D_{0}^{-1}\left(\sum^{n}_{k = 1} \phi^{*}_{k}F_{k}\right)\left(\sum^{n}_{k = 1} F^{*}_{k}\phi_{k}\right)D_{0}^{-1}\psi^{*}_{n}\right] &  \\
& & \\
= & \tr \left[ D_{0}^{-1}\left(\sum^{n}_{k = 1} \phi^{*}_{k}F_{k}\right)\left(\sum^{n}_{k = 1} F^{*}_{k}\phi_{k}\right)D_{0}^{-1}\psi^{*}_{n}\psi_{n}\right] &  \\
\end{array}
$$
and then, from Cauchy-Schwarz Inequality applied to the Frobenius inner product, that
$$
\begin{array}{lcl}
& \tr \left[\left(\sum^{n}_{k = 1} F^{*}_{k}\phi_{k}\right)D_{0}^{-1}\psi^{*}_{n}\psi_{n} D_{0}^{-1}\left(\sum^{n}_{k = 1} \phi^{*}_{k}F_{k}\right)\right] &  \\ 
& & \\
 \leq &  \tr \left[ D_{0}^{-1}\left(\sum^{n}_{k = 1} \phi^{*}_{k}F_{k}\right)\left(\sum^{n}_{k = 1} F^{*}_{k}\phi_{k}\right)D_{0}^{-1} \right] \tr \left[ \psi^{*}_{n}\psi_{n}\right] & \\
& & \\
\leq & \tr \left[ \left(\sum^{n}_{k = 1} \phi^{*}_{k}F_{k}\right)\left(\sum^{n}_{k = 1} F^{*}_{k}\phi_{k}\right)\right] \tr \left[ D_{0}^{-2} \right] \tr \left[ \psi^{*}_{n}\psi_{n}\right].&
\end{array}
$$
Hence, for each $L\ge 1$,
$$
\begin{array}{lll}
\left\| \psi\alpha \right\|^{2}_{L} &\leq &  
\sum^{\left\lfloor  L \right\rfloor}_{n = 1} \left\|\psi_{n}\right\|_F^{2} \left\|D_{0}^{-1}\right\|_F^{2} \tr \left[\left(\sum^{n}_{k = 1} F^{*}_{k}\phi_{k}\right)\left(\sum^{n}_{k = 1} \phi^{*}_{k}F_{k}\right)\right]  \\
& & \\
&  + & \left\|\psi_{\left\lfloor  L \right\rfloor + 1}\right\|_F^{2} \left\|D_{0}^{-1}\right\|^{2}_F ( L - \left\lfloor  L \right\rfloor ) \tr \left[\left(\sum^{\left\lfloor  L \right\rfloor + 1}_{k = 1} F^{*}_{k}\phi_{k}\right)\left(\sum^{\left\lfloor  L \right\rfloor + 1}_{k = 1} \phi^{*}_{k}F_{k}\right)\right] \\
& & \\
& = & \left\|D_{0}^{-1}\right\|_F^{2} \sum^{\left\lfloor  L \right\rfloor}_{n = 1} \left\|\psi_{n}\right\|_F^{2} \left(\sum^{n}_{p = 1} \sum^{n}_{q = 1}  \tr \left[F^{*}_{p}\phi_{p}\phi^{*}_{q}F_{q}\right]\right)\\
& & \\
&  + & \left\|\psi_{\left\lfloor  L \right\rfloor + 1}\right\|_F^{2} \left\|D_{0}^{-1}\right\|_F^{2} ( L - \left\lfloor  L \right\rfloor ) \left(\sum^{\left\lfloor  L \right\rfloor + 1}_{p = 1} \sum^{\left\lfloor  L \right\rfloor + 1}_{q = 1}  \tr \left[F^{*}_{p}\phi_{p}\phi^{*}_{q}F_{q}\right]\right).
\end{array}
$$

Again, by Cauchy-Schwarz Inequality (applied to the Frobenius inner product, in the first inequality, and then to the inner product of $l^2(\N;\R)$, in the second one), one has for each $n \in \mathbb{N}$ that
$$
\begin{array}{lll}
  & \sum^{n}_{p = 1} \sum^{n}_{q = 1}  \tr \left[F^{*}_{p}\phi_{p}\phi^{*}_{q}F_{q}\right] \le &\\
&&\\  
& \sum^{n}_{p = 1} \sum^{n}_{q = 1}  \sqrt{\tr \left[F^{*}_{p}F_{p}\right] \tr \left[\phi_{p}\phi^{*}_{p}\right]  \tr \left[ \phi^{*}_{q}\phi_{q}\right] \tr \left[F_{q}F^{*}_{q}\right]} & \\
& & \\
& \leq  \left(\sum^{n}_{p = 1} \sqrt{\tr \left[F^{*}_{p}F_{p}\right] \tr \left[\phi_{p}\phi^{*}_{p}\right] } \right)\left(\sum^{n}_{q = 1} \sqrt{\tr \left[ \phi^{*}_{q}\phi_{q}\right] \tr \left[F_{q}F^{*}_{q}\right]}\right)&\\
& & \\
& = \left(\left\|F\right\|_{n} \right)^{2}  \left(\left\|\phi\right\|_{n} \right)^{2},&
\end{array}
$$
and consequently, that 
$$
\begin{array}{lll}
\left\| \psi\alpha \right\|^{2}_{L} & \leq &  \left\|D_{0}^{-1}\right\|^{2}_F \sum^{\left\lfloor  L \right\rfloor}_{n = 1} \left\|\psi_{n}\right\|^{2}_F \left(\left\|F\right\|_{n} \right)^{2}  \left(\left\|\phi\right\|_{n} \right)^{2} \\
& & \\
& & + \left\|\psi_{\left\lfloor  L \right\rfloor + 1}\right\|^{2}_F \left\|D_{0}^{-1}\right\|^{2}_F ( L - \left\lfloor  L \right\rfloor ) \left(\left\|F\right\|_{\left\lfloor  L \right\rfloor + 1} \right)^{2}  \left(\left\|\phi\right\|_{\left\lfloor  L \right\rfloor + 1} \right)^{2} \\
& & \\
& \leq & \left\|D_{0}^{-1}\right\|^{2}_F \left(\left\|\psi\right\|_{\left\lfloor L \right\rfloor}\right)^{2} \left(\left\|F\right\|_{\left\lfloor  L \right\rfloor + 1} \right)^{2}  \left(\left\|\phi\right\|_{\left\lfloor  L \right\rfloor} \right)^{2} \\
& & \\
& & + \left(\left\|\psi\right\|_{\left\lfloor  L \right\rfloor + 1}\right)^{2} \left\|D_{0}^{-1}\right\|^{2}_F ( L - \left\lfloor  L \right\rfloor ) \left(\left\|F\right\|_{\left\lfloor  L \right\rfloor + 1} \right)^{2}  \left(\left\|\phi\right\|_{\left\lfloor  L \right\rfloor + 1} \right)^{2}.
\end{array}
$$
 Thus, if one lets $L$ be as in \eqref{def.l}, one gets
$$
\begin{array}{lll}
\left\| \psi\alpha \right\|^{2}_{L} & \leq & 
\left\|D_{0}^{-1}\right\|^{2}_F \left\|F\right\|_{\left\lfloor  L \right\rfloor + 1}^{2}  \left[\left\|\phi\right\|_{\left\lfloor  L \right\rfloor}^2\left\|\psi\right\|_{\left\lfloor  L \right\rfloor}^{2} +
  ( L - \left\lfloor  L \right\rfloor )\left\|\phi\right\|_{\left\lfloor  L \right\rfloor + 1}^{2} \left\|\psi\right\|_{\left\lfloor  L \right\rfloor + 1}^{2}\right] \\
& & \\
& \leq & \dfrac{\left\|F\right\|_{\left\lfloor  L \right\rfloor + 1}^2}{y^2}.
\end{array}
$$

In the same fashion, one can prove that 
$$
\left\| \phi\beta \right\|_{L} \leq \frac{\left\| F \right\|_{\left\lfloor  L \right\rfloor + 1}}{y}.
$$

The next step consists in estimating $\Vert F\Vert_L$ from below. 
It follows from \eqref{eq.f.alpha} and from Triangular Inequality that 
$$
\begin{array}{lll}
\left\| F \right\|_{L} & \geq & \left\| \psi - \phi M^{\phi}D_{0} \right\|_{L} - \left\| iy\psi \alpha \right\|_{L} - \left\|iy \phi \beta \right\|_{L} \\
&& \\
&\geq &\left\| \psi - \phi M^{\phi}D_{0} \right\|_{L} - 2 y \dfrac{\left\| F \right\|_{\left\lfloor L \right\rfloor + 1}}{y},
\end{array}
$$
so 
$$
\left\| \psi - \phi M^{\phi}D_{0} \right\|_{L}  
\leq  3 \left\| F \right\|_{\left\lfloor L \right\rfloor + 1}.
$$

By squaring both sides of the previous inequality, one gets

\begin{equation}
\label{desi.para.m.f}
\begin{array}{lll}
 9  \left\| F \right\|_{\left\lfloor L \right\rfloor + 1}^{2} &\geq &
\left\| \psi \right\|^{2}_{L} + \left\| \phi M^{\phi}D_{0} \right\|^{2}_{L} - \left(L - \left\lfloor L \right\rfloor \right) \tr \left[ \psi^{*}_{\left\lfloor L \right\rfloor + 1}\phi_{\left\lfloor L \right\rfloor + 1}M^{\phi}D_{0} \right] \\
& & \\
&  - &\sum^{\left\lfloor L \right\rfloor}_{n = 1}  \tr \left[ \psi^{*}_{n}\phi_{n}M^{\phi}D_{0} \right] - \sum^{\left\lfloor L \right\rfloor}_{n = 1}  \tr \left[ D_{0}(M^{\phi})^{*}\phi^{*}_{n}\psi_{n}  \right] \\
& & \\
&   - & \left(L - \left\lfloor L \right\rfloor \right) \tr \left[ D_{0}(M^{\phi})^{*}\phi^{*}_{\left\lfloor L \right\rfloor + 1}\psi_{\left\lfloor L \right\rfloor + 1}  \right].
\end{array}
\end{equation}

Now, it follows from Cauchy-Schwarz Inequality (applied to the Frobenius inner product) that, for each $n \in \mathbb{N}$,
$$
\begin{array}{lll}
& \tr \left[ D_{0}(M^{\phi})^{*}\phi^{*}_{n}\psi_{n}  \right] & 
\\
\leq & \sqrt{\tr \left[ D^{*}_{0}D_{0} \right]}\sqrt{\tr \left[ \psi^{*}_{n}\phi_{n}M^{\phi}  (M^{\phi})^{*}\phi^{*}_{n}\psi_{n} \right]}& \\
& & \\
 = & \left\| D_{0} \right\|_F \sqrt{\tr \left[ \psi_{n}\psi^{*}_{n} \phi_{n}M^{\phi}  (M^{\phi})^{*}\phi^{*}_{n} \right]}& \\
& & \\
 \leq & \left\| D_{0} \right\|_F \sqrt{ \sqrt{\tr \left[ \psi^{*}_{n}\psi_{n} \psi_{n}\psi^{*}_{n} \right]} \sqrt{ \tr \left[   \phi_{n} M^{\phi} (M^{\phi})^{*}\phi^{*}_{n} \phi_{n}M^{\phi}  (M^{\phi})^{*}\phi^{*}_{n}\right]}}&\\
& & \\
 \leq & \left\| D_{0} \right\|_F \sqrt{ \tr \left[ \psi^{*}_{n} \right] \tr \left[   \phi_{n} M^{\phi} (M^{\phi})^{*}\phi^{*}_{n} \right]}  &\\
& & \\
 = & \left\| D_{0} \right\|_F \left\| \psi_{n} \right\|_F \sqrt{ \tr \left[  \phi^{*}_{n} \phi_{n} M^{\phi} (M^{\phi})^{*} \right]} &\\
& & \\
 \leq & \left\| D_{0} \right\|_F \left\| \psi_{n} \right\|_F \sqrt{ \sqrt{\tr \left[  \phi^{*}_{n} \phi_{n} \phi_{n}\phi^{*}_{n} \right]} \sqrt{\tr \left[  M^{\phi} (M^{\phi})^{*}(M^{\phi})^{*} M^{\phi}  \right]}}  &\\
& & \\
 \leq & \left\| D_{0} \right\|_F \left\| \psi_{n} \right\|_F \sqrt{ \tr \left[  \phi^{*}_{n} \phi_{n} \right] \tr \left[  M^{\phi} (M^{\phi})^{*}  \right]}  & \\
& & \\
= &  \left\| D_{0} \right\|_F \left\| \psi_{n} \right\|_F  \left\|  \phi_{n} \right\|_F \left\| M^{\phi} \right\|_F,&
\end{array}
$$
and so, by the same reasoning, that
\[\tr \left[ \psi^{*}_{n}\phi_{n}M^{\phi}D_{0} \right]\le \left\| D_{0} \right\|_F \left\| \psi_{n} \right\|_F  \left\|  \phi_{n} \right\|_F \Vert M^{\phi} \Vert_F.\]

Then, by letting $L$ be as in~\eqref{def.l}, one has 
$$
\begin{array}{lll}
& \sum^{\left\lfloor L \right\rfloor}_{n = 1}  \tr \left[ \psi^{*}_{n}\phi_{n}M^\phi D_{0} \right] + \left(L - \left\lfloor L \right\rfloor \right) \tr \left[ \psi^{*}_{\left\lfloor L \right\rfloor + 1}\phi_{\left\lfloor L \right\rfloor + 1}M^{\phi}D_{0} \right] &   \\
& & \\
\leq & \left\| D_{0} \right\|_F \left\| M^{\phi} \right\|_F \left( \left\| \psi \right\|_{\left\lfloor L \right\rfloor} \left\|  \phi \right\|_{\left\lfloor L \right\rfloor} + \left(L - \left\lfloor L \right\rfloor \right)\left\| \psi_{\left\lfloor L \right\rfloor + 1} \right\| \left\|  \phi_{\left\lfloor L\right\rfloor+1} \right\|  \right) &  \\
& & \\
\leq &  2\left\| D_{0} \right\|_F\left\| M^{\phi} \right\|_F\left\| \psi \right\|_{\left\lfloor L \right\rfloor + 1} \left\|  \phi \right\|_{\left\lfloor L \right\rfloor+1} = \left\| D_{0} \right\|_F\left\| M^{\phi} \right\|_F \dfrac{1}{y \left\| D^{-1}_{0} \right\|_F}.&  
\end{array}
$$
Now, by applying such estimates to the right-hand side of relation~\eqref{desi.para.m.f}, it follows that
\begin{equation}
\label{des.f.m.d}
9  \left\| F \right\|_{\left\lfloor L \right\rfloor + 1}^{2}  \geq  \left\| \psi \right\|^{2}_{L} + \left\| \phi M^{\phi}D_{0} \right\|^{2}_{L} - 2 \left\| D_{0} \right\|_F\left\| M^\phi \right\|_F \frac{1}{y \left\| D^{-1}_{0} \right\|_F}.
\end{equation}

Therefore, by noting that for each pair of matrices $A,B\in M(l;\C)$, one has from Definition~\ref{frobeniusnorm} and Proposition~\ref{traceFrobenius} that $s_1[A^\ast A]\le\Vert A\Vert^2_F$, $\tr[A^\ast A]=\sum_{k=1}^ls_k[A^\ast A]$, and from Proposition~\ref{Th6.7Haeal}-c that $s_l[A]s_1[B]\le s_1[AB]$, it follows that
$$
\begin{array}{lll}
\left\| \phi M^{\phi}D_{0} \right\|^{2}_{L} & = & \sum^{\left\lfloor L \right\rfloor}_{n = 1} \tr\left[ \phi_{n} M^{\phi}D_{0} D_{0} (M^{\phi})^{*}\phi_{n}^{*} \right] \\
& & \\
& & + \left(L - \left\lfloor L \right\rfloor \right) \tr\left[ \phi_{\left\lfloor L \right\rfloor + 1} M^{\phi}D_{0} D_{0} (M^{\phi})^{*}\phi_{\left\lfloor L \right\rfloor + 1}^{*} \right]\\
& & \\
& \geq & l^{-1}\left\|M^{\phi}\right\|^{2}_F s_{l}[D^{2}_{0}] \sum^{\left\lfloor L \right\rfloor}_{n = 1}  s_{l}\left[ \phi_{n}\phi_{n}^{*} \right] \\
&& \\
& & + l^{-1}\left(L - \left\lfloor L \right\rfloor \right) \left\|M^{\phi}\right\|^{2}_F s_{l}[D^{2}_{0}] s_{l}\left[ \phi_{\left\lfloor L \right\rfloor + 1} \phi_{\left\lfloor L \right\rfloor + 1}^{*} \right]\\
& & \\
& = & l^{-1}\left\|M^{\phi}\right\|^{2}_F s_{l}[D^{2}_{0}] s_{l}[\phi]^{2}_{L};
\end{array}
$$
thus, by relation \eqref{des.f.m.d}, one gets 
\begin{equation}
\label{des.m.f.2}
9 \left\| F \right\|_{\left\lfloor L \right\rfloor + 1}^{2}  \geq  \left\| \psi \right\|^{2}_{L}  - 2 \left\| D_{0} \right\|_F\left\| M^{\phi} \right\|_F \frac{1}{y \left\| D^{-1}_{0} \right\|_F} + l^{-1}\left\|M^{\phi}\right\|^{2}_F s_{l}[D^{2}_{0}] s_{l}[\phi]^{2}_{L}.
\end{equation}

The next step consists in bounding $\Vert F\Vert_L$ from above. It follows from Proposition \ref{prop.m}-$b$ that
\begin{eqnarray*}
\frac{\tr[D_{0} \Im[M^{\phi}] D_{0}]}{y} & \geq  &  \left\| F \right\|_{\left\lfloor L \right\rfloor + 1}^{2}.
\end{eqnarray*}
By Cauchy-Schwarz Inequality, one has
$$
\tr[D_{0} \Im[M^{\phi}] D_{0}] = \tr[D^{2}_{0} \Im[M^{\phi}]] \leq \left\| D^{2}_{0} \right\|_F \left\| \Im[M^{\phi}] \right\|_F,
$$
and since $\Im[M^{\phi}] = \frac{1}{2i} (M^{\phi} - \overline{M^{\phi}})$, it follows that $\left\| \Im[M^{\phi}] \right\|_F \leq \left\| M^{\phi} \right\|_F$. Thus,
\begin{equation}\label{est.F.above}
\frac{\Vert D_{0}^2\Vert_F\Vert M^{\phi}\Vert_F}{y} \geq   \left\| F \right\|_{\left\lfloor L \right\rfloor + 1}^{2}.
\end{equation}

At last, one gets from relations \eqref{des.m.f.2} and \eqref{est.F.above} a quadratic inequality for $\left\| M^{\phi} \right\|$:
$$
\frac{9}{y} \left\| D^{2}_{0} \right\|_F \left\| M^{\phi} \right\|_F  \geq  \left\| \psi \right\|^{2}_{L}  -  2 \left\| D_{0} \right\|_F\left\| M^{\phi} \right\|_F \frac{1}{y \left\| D^{-1}_{0} \right\|_F} + l^{-1}\left\|M^{\phi}\right\|^{2}_F s_{l}[D^{2}_{0}] s_{l}[\phi]^{2}_{L},
$$
that is, 
$$
y \left\| \psi \right\|^{2}_{L}  - \left(2 \frac{ \left\| D_{0} \right\|_F}{\left\| D^{-1}_{0} \right\|_F} + 9 \left\| D^{2}_{0} \right\|_F  \right) \left\| M^{\phi} \right\|_F + yl^{-1} \left\|M^{\phi}\right\|_F^{2} s_{l}[D^{2}_{0}] s_{l}[\phi]^{2}_{L} \leq 0,
$$
which is an inequality of the form
$$
a \left\|M^{\phi}\right\|^{2}_F + b \left\|M^{\phi}\right\|_F + c \leq 0,
$$
with
$$
\left\{
\begin{array}{l}
a  :=  yl^{-1}  \left(s_{l}[D^{2}_{0}]s_{l}[\phi]^{2}_{L}\right),  \\
\\
b  :=  - \left( 2 \dfrac{\left\| D_{0}\right\|_F}{\left\|D^{-1}_{0}\right\|_F}  +  9\left\| D_{0} \right\|^{2}_F \right), \\
\\
c  :=  y \left\| \psi \right\|^{2}_{L}.
\end{array}\right.
$$

By computing its discriminant, it follows from relation~\eqref{def.l} that
$$
\begin{array}{lll}
\Delta & = & \left( 2 \dfrac{\left\| D_{0}\right\|_F}{\left\|D^{-1}_{0}\right\|_F}  +  9\left\| D_{0} \right\|^{2}_F \right)^{2} - 4 y^{2}l^{-2} \left(s_{l}[D^{2}_{0}]\right) \left(s_{l}[\phi]^{2}_{L}\right) \left\| \psi \right\|^{2}_{L} \\
& \geq & \left( 2 \dfrac{\left\| D_{0}\right\|_F}{\left\|D^{-1}_{0}\right\|_F}  +  9\left\| D_{0} \right\|^{2}_F \right)^{2} - \left( \dfrac{s_{l}[D^{2}_{0}]}{\left\|D^{-1}_{0}\right\|^{2}}  \right) \\
& > & 0, 
\end{array}
$$
and so, the correspondent quadratic equation has two distinct solutions. Finally, one gets
$$
 \frac{-c}{b} < \left\| M^{\phi} \right\|_F < \frac{-2 b}{a}, 
$$
that is,
$$
\left( - \frac{1}{b} \right)  \left( \frac{1}{\left\|D^{-1}_{0}\right\|_F} \right) \left( \frac{\left\| \psi \right\|_{L}}{\left\| \phi \right\|_{L}}\right)     < \left\| M^{\phi} \right\|_F < - 2l b \left( \frac{\left\|D^{-1}_{0}\right\|_F}{s_{l}[D^{2}_{0}]}  \right) \left( \frac{\left\| \phi \right\|^{2}_{L}}{ s_{l}[\phi]^{2}_{L} }  \right) \left( \frac{\left\| \psi \right\|_{L}}{\left\| \phi \right\|_{L}} \right).
$$

The result follows now by setting 
$$
\left\{
\begin{array}{lll}
k_{1} & := & - \left( \dfrac{C_1^2}{C_2 b} \right)  \left( \dfrac{1}{\left\|D^{-1}_{0}\right\|_F} \right), \\
& & \\ 
k_{2} & := & - 2l b \left( \dfrac{C_2\left\|D^{-1}_{0}\right\|_F}{C_1^2s_{l}[D^{2}_{0}]}  \right), 
\end{array}\right.
$$
where $C_2\ge C_1>0$ are such that, for each $A\in M(l,\C)$,
\[C_1\Vert A\Vert_F\le \Vert A\Vert\le C_2\Vert A\Vert_F\]
(recall that any two norms in $M(l,\C)$ are equivalent).
\end{proof1}

\section{Criteria for the Absolutely Continuous Spectral Components}
\label{crit.abs.cont.sp.comp.}
 
From now on, we assume that the singular values of the matrices $D_{n}$, $n\in\Z$, satisfy relation~\eqref{eq.limit.dn.1}, that is,
\[0 < \inf_{n \in \Z} s_{l}[D_{n}] \leq \sup_{n \in \Z} s_{1}[D_{n}] < \infty. 
\]

Naturally, it follows from the first condition in~\eqref{eq.limit.dn.1} that these matrices are non-singular (in this case, the operator $H$ is called non-singular; see \cite{marx15}), that is,  $\det(D_{n}) \neq 0$ for every $n \in \Z$. Note, however, that we are assuming a stronger condition, namely, that the bilateral sequence $(s_{l}[D_{n}])_{n \in \mathbb{Z}}$ is uniformly bounded from below by a positive number.


\subsection{Characterization of the Absolutely Continuous Spectrum}
\label{sec.mat.tra}

Our main goal in this section is to prove that, for each $r\in\{1,\ldots,l - 1\}$ (respectively, $r = l$), the essential closure (with respect to the Lebesgue measure) of the set $\mathcal{S}_{r} \setminus \mathcal{S}_{r + 1}$ (respectively, $\mathcal{S}_{l}$), where $\mathcal{S}_{r}$ is given by \eqref{def.val.r.matriz}, is a minimal support to the absolutely continuous component of multiplicity $r$ (respectively, $l$) of the operator $H^{\phi}_+$. Moreover, since $\sigma_{ac, r}(H^{\phi}_+) = \sigma_{ac, r}(H^{\psi}_+)$ (by Proposition~\ref{prop.supor.pertur}), $\overline{\mathcal{S}_{r} \setminus \mathcal{S}_{r + 1}}^{ess}$ is a minimal support to $\sigma_{ac, r}(H^{\psi}_+)$. In fact, the same can be said about any self-adjoint extension of $H^{+}_{min}$ (which is particularly true for the extensions given by \eqref{def.extensoes1} and \eqref{def.extensoes2}), since they can be seen as finite rank perturbations of $H^\phi_+$.

\begin{remark} Although we just present the results for the absolutely continuous spectrum of $H^\phi_+$, the same conclusions are valid for the absolutely continuous spectrum of $H^\phi_-$.
\end{remark}

The first step in our analysis consists in showing that if $z\in\rho(H^\phi_+)$, then there exist a solution to the eigenvalue equation~\eqref{eq.aut.Jacobi} for $H^\phi_+$ that decays exponentially fast as $n\rightarrow\infty$, and a solution to the same equation that grows exponentially fast as $n\rightarrow\infty$ (here, we follow an argument presented in Section~2.2 in~\cite{teschl00}). 

\begin{3}
\label{lema.dcai.exp.ess}
Let $H^{\phi}_+$ be given by \eqref{def.oper.diri} and let, for each $k\in\{1,\ldots, l\}$, $H^{\phi}_{k}$ be the restriction of $H^{\phi}$ to the subspace spanned by $\textbf{e}_{1, k}$. 
Then, for each $z \in \mathbb{C} \setminus \sigma(H^{\phi}_{k})$, there exists a solution to the eigenvalue equation for $H^{\phi}_{k}$ at $z$ that goes to zero exponentially fast as $n\rightarrow\infty$.
\end{3}
\begin{proof}
Let $k\in\{1,\ldots, l\}$ and $z \in \mathbb{C} \setminus \sigma(H^{\phi}_{k})$. We will show that for each $\alpha>0$, there exists a constant $C>0$ such that, for each $n \in \mathbb{N}$,
$$
\left| g(1, n) e^{\alpha n} \right|\le C,
$$ 
with 
$$
g(1, n) := \left\langle  (\textbf{e}_{1, k})_{1},  ( (H^{\phi}_{k} - z)^{-1} \textbf{e}_{1, k})_{n}  \right\rangle_{\mathbb{C}^{l}}
$$
(that is, the Green Function of $H^{\phi}_{k}$ at $(1,n)$).
 
So, let $\alpha > 0$ and consider the operators $P_\alpha,Q_\alpha:(\mathbb{C}^{l})^{\mathbb{N}}\rightarrow(\mathbb{C}^{l})^{\mathbb{N}}$ given by 
$$
\begin{array}{lll}
(Q_{\alpha}\textbf{u})_{n} & := & (e^{\alpha} - 1)D_{n}\textbf{u}_{n + 1} + (e^{-\alpha} - 1)D_{n - 1}\textbf{u}_{n - 1},\\
& & \\
(P_{\alpha}\textbf{u})_{n} & := & e^{\alpha n}\textbf{u}_{n}.
\end{array}
$$
Note that, for each $n \in \mathbb{N}$,
$$
e^{\alpha (n - 1)}g(1, n) = \left\langle  (\textbf{e}_{1, k})_{1},  (P_{-\alpha} (H^{\phi}_{k} - z)^{-1} P_{\alpha} \textbf{e}_{1, k})_{n}  \right\rangle_{\mathbb{C}^{l}}.
$$

On the other hand, for each $\textbf{u}\in\dom(H^\phi_k)$, one has
$$
\left(P_{-\alpha} (H^{\phi}_{k} - z) P_{\alpha}\right)\textbf{u} = \left(H^{\phi}_{k} - z + Q_{\alpha}\right)\textbf{u},
$$
from which follows that 
$$
e^{\alpha (n - 1)}g(1, n) = \left\langle  (\textbf{e}_{1, k},  ((H^{\phi}_{k} - z + Q_{\alpha})^{-1} \textbf{e}_{1, k})_{n}  \right\rangle_{\mathbb{C}^{l}}.
$$

Hence, if one proves that $(H^{\phi}_{k} - z + Q_{\alpha})^{-1}$ is a bounded operator, one concludes that the sequence $\left(e^{\alpha (n - 1)}g(1, n)\right)_{n \in \mathbb{N}}$ is also bounded. Firstly, note that 
$$
\left\| Q_{\alpha} \right\| \leq 2 (e^{\alpha} - 1) \sup_{n} \left\| D_{n} \right\|.
$$

Now, set $\delta := \dist(z, \sigma(H^{\phi}_{k}))$ and let $\epsilon > 0$ be such that
$$
\alpha = \log \left(1 + \frac{(1 - \epsilon)}{2\sup_{n} \left\| D_{n} \right\|} \right), 
$$
from which follows that $\left\| Q_{\alpha} \right\| \leq (1 - \epsilon)\delta$. The second resolvent identity establishes that
$$
(H^{\phi}_{k} - z + Q_{\alpha})^{-1} = (H^{\phi}_{k} - z)^{-1} + (H^{\phi}_{k} - z + Q_{\alpha})^{-1}(- Q_{\alpha})(H^{\phi}_{k} - z)^{-1};
$$
therefore, 
$$
\left\|(H^{\phi}_{k} - z + Q_{\alpha})^{-1}\right\| \leq \left\| (H^{\phi}_{k} - z)^{-1} \right\|  +  \left\| (H^{\phi}_{k} - z)^{-1} \right\| \left\| Q_{\alpha} \right\|   \left\| (H^{\phi}_{k} - z + Q_{\alpha})^{-1} \right\|, 
$$
and, given that $\left\| (H^{\phi}_{k} - z)^{-1} \right\| = \delta^{-1}$, one gets 
$$
\left\|(H^{\phi}_{k} - z + Q_{\alpha})^{-1}\right\| \leq \frac{1}{\delta\epsilon}.
$$

Finally, by letting $C:=e^\alpha/(\delta\epsilon)$, it follows that for each $n\in\N$, $|e^{\alpha n}g(1, n)|\le C$, and we are done. 
\end{proof}

\begin{3}
\label{lema.sol.cresc.psi}
Let $H^{\phi}_+$ and $H^{\psi}_+$ be the operators given by \eqref{def.oper.diri} and \eqref{def.oper.neu}, respectively, and let, for each $k\in\{1,\ldots,l\}$,  $H^{\phi,\psi}_{k}$ be as in the statement of Lemma \ref{lema.dcai.exp.ess}. If $z \in \mathbb{C} \setminus \sigma(H^{\phi}_{k})$, then there exists a solution to the eigenvalue equation for $H^{\psi}_k$ at $z$ which grows exponentially fast as $n\rightarrow\infty$.
\end{3}
\begin{proof}
  Let $k\in\{1,\ldots, l\}$ and $z \in \mathbb{C} \setminus \sigma(H^{\phi}_{k})$. It follows from Lemma \ref{lema.dcai.exp.ess} that there exists a solution $\textbf{u}$ to the eigenvalue equation for $H_{k}^{\phi}$ at $z$ which goes to zero exponentially fast as $n\rightarrow\infty$.

  Now, if $\textbf{v}$ is any solution to the eigenvalue equation for $H^{\psi}_k$ at $z$, then necessarily $\textbf{v}_{1} = \textbf{0}$, and so 
$$
W_{[\textbf{u}, \textbf{v}]}(1) =  \left\langle D_{0} \textbf{u}_{1}, \bar{\textbf{v}}_{0} \right\rangle_{\mathbb{C}^{l}} - \left\langle D_{0} \textbf{v}_{1}, \overline{\textbf{u}}_{0}  \right\rangle_{\mathbb{C}^{l}} = \textbf{u}^{t}_{1}D_{0}\textbf{v}_{0}.
$$

Since $D_{0}$ is invertible, there exists $\textbf{v}$ so that $W_{[\textbf{u}, \textbf{v}]}(1) = C \neq 0$. It follows now from the constancy of the Wronskian 
that for each $n \in \N$, 
$$
W_{[\textbf{u}, \textbf{v}]}(n) =  \left\langle D_{n - 1}\textbf{u}_{n}, \bar{\textbf{v}}_{n - 1} \right\rangle_{\mathbb{C}^{l}} - \left\langle D_{n - 1} \textbf{v}_{n}, \overline{\textbf{u}}_{n - 1}  \right\rangle_{\mathbb{C}^{l}} = C
$$
(in order to obtain the result, one has to use, in Green Formula, that $\textbf{u}_0=0$ and that, for each $j\in\{1,\ldots,l\}$, $H_j^\phi\textbf{v}=H_j^\psi\textbf{v}-D_0\langle\textbf{v},\textbf{e}_{1,j}\rangle_{\mathbb{C}^{l}}\textbf{e}_{0,j}$), and then, by Cauchy-Schwarz Inequality, that
$$
\begin{array}{lll}
\left| C \right| & \leq & \left\| D_{n - 1} \right\| \left( \left\| \textbf{u}_{n} \right\| \left\| \bar{\textbf{v}}_{n - 1} \right\| + \left\| \textbf{v}_{n} \right\| \left\| \overline{\textbf{u}}_{n - 1} \right\| \right)\\
& & \\
& \leq & \left(\sup_{n \in \N} \left\| D_{n - 1} \right\|\right) \left( \left\| \textbf{u}_{n} \right\| \left\| \bar{\textbf{v}}_{n - 1} \right\| + \left\| \textbf{v}_{n} \right\| \left\| \overline{\textbf{u}}_{n - 1} \right\| \right).
\end{array}
$$

Therefore, $\left\| \textbf{v}_{n} \right\|$ grows exponentially fast as $n\rightarrow\infty$.

\end{proof}

\begin{7}
\label{prop.suporte.cesaro}
Let $r\in\{1,\ldots,l\}$. Then, the restriction of the absolutely continuous spectrum of the operator $H^{\phi}_+$, given by \eqref{def.oper.diri}, to the set $\mathcal{S}_{r}$, given by \eqref{def.val.r.matriz}, has multiplicity at least $r$.
\end{7}
\begin{proof}
Let $x \in \mathbb{R}$, consider the spectral basis $\{\textbf{e}_{1, 1}, \textbf{e}_{1, 2}, \ldots, \textbf{e}_{1, l}\}$ for $H^{\phi}_+$ and let, for each $j\in\{1,\ldots,l\}$, $H^{\phi}_{j}$ be the restriction of $H^{\phi}_+$ to the subspace spanned by $\textbf{e}_{1, j}$. Set
$$
\mathcal{T}:=\bigcup_{(k_{1}, \ldots, k_{r})}  \bigcap_{j = 1}^{r} \sigma(H^{\phi}_{k_{j}}),
$$
where $(k_{1}, \ldots, k_{r})$ is a $r$-uple of elements of $\{1,2,\ldots,l\}$ such that for each $i\neq j\in\{1,\ldots,r\}$, $k_i\neq k_j$, and note that $\overline{\mathcal{T}}^{ess}$ is a minimal support of the absolutely continuous spectrum of $H^\phi_+$ of multiplicity at least $r$. Suppose now that $x\notin\mathcal{T}$; then, for each $r$-uple $(k_{1}, \ldots, k_{r})$ as defined above, there exists $k_{i}$ such that 
$$
x \notin \sigma(H^{\phi}_{k_{i}}).
$$

If one takes, for instance, the $r$-uple $(1, 2, \ldots, r)$, there exists $k_{1} \in \{1, 2, \ldots, r\}$ such that
$$
x \notin \sigma(H^{\phi}_{k_{1}}).
$$
Consider now the $r$-uple obtained by replacing $k_{1}$ by 
any $p_1\in\{r + 1,r+2,\ldots,l\}$, that is, the $r$-uple $(1, 2, \ldots, k_{1} - 1, p_1, k_{1} + 1, \ldots, r - 1,  r)$. 
By applying the hypothesis to this new $r$-uple, it follows that there exists $k_{2} \in \{1, 2, \ldots, k_{1} - 1, p_1, k_{1} + 1, \ldots, r - 1,  r \}$  such that
$$
x \notin\sigma(H^{\phi}_{k_{2}})
$$
(note also that $k_{2} \neq k_{1}$). One may now consider the $r$-uple obtained from the previous one by replacing $k_{2}$ by any $p_2\in\{r + 1,r+2,\ldots,l\}\setminus\{p_1\}$.

In this way, by repeating the process, one concludes that there are at least $l - r + 1$  distinct values of $k_{i}$ such that
$$
x \notin\sigma(H^{\phi}_{k_{i}}).
$$

Thus, $x\in\bigcup_{(k_1,\ldots,k_r)}\bigcap_{i=1}^{l-r+1}(\C\setminus\sigma(H^{\phi}_{k_{i}}))$. It follows from Lemma \ref{lema.sol.cresc.psi} and from the fact that $H^\psi_+$ is in the limit point case at $+\infty$ that there exist $l - r + 1$ distinct solutions to the eigenvalue equation for $H^{\psi}$ at $x$ which grow exponentially fast as $n\rightarrow\infty$. Consequently, given that the columns of $(\psi_{n}(x))$ are solutions to the eigenvalue equation~\eqref{eq.aut.Jacobi} that satisfy the Neumann boundary condition, it follows from 
Theorem $6.7$ in \cite{petz2014} that
$$
\lim_{n \rightarrow \infty} s_{l - r + 1}[\psi_{n}(x)] = \infty,
$$
and so, that 
$$
\liminf_{L \rightarrow \infty} \frac{1}{L} \sum^{L}_{n = 1} s^{2}_{l - r + 1}[\psi_{n}(x)] = \infty.
$$

One concludes, therefore, that $\overline{\mathcal{S}_r}^{ess}\subset\overline{\mathcal{T}}^{ess}$.
\end{proof}


\

We proceed now to the proof of the converse of Theorem \ref{prop.suporte.cesaro}.
 
\begin{3}
\label{lema.g.ident}
Let  $G^{\phi}$ be the matrix-valued Green Function of $H^{\phi}_+$, given by Proposition \ref{porp.def.green}. Then, for each $n \in \N$,
$$
\lim_{y \downarrow 0} \Im G^{\phi}(n, n; x + iy) d\kappa(x) = \lim_{y \downarrow 0} 
\phi_{n}(x)\Im\left[M^{\phi}_+(x + iy)\right]\phi_{n}^{t}(x) d\kappa(x),
$$ 
where $M^{\phi}_+$ is the matrix-valued Weyl-Titchmarsh function given by Definition \ref{def.m.weyl1}. 
\end{3}
\begin{proof}
We omit the dependence on $(+)$ throughout the proof. By definition, the Green Function of the operator $H^{\phi}$, for each $n \in \mathbb{N}$, is given by
$$
G(n, n; x + iy) = - \phi_{n} (x + iy) D_{0}^{-1} F^{t}_{n}(x + iy). 
$$
Since one has from relation \eqref{equa.m.jost} and Proposition~\ref{prop.m}-a) that for each $n\in\N$, $F^{t}_{n}(x + iy) = \psi^{t}_{n}(x + iy) - D_{0}M^{\phi}(x + iy)\phi_{n}^{t}(x + iy)$, it follows that 
$$
\begin{array}{lll}
G(n, n; x + iy) & = & - \phi_{n}(x + iy)D_{0}^{-1} (\psi^{t}_{n}(x + iy) - D_{0}M^{\phi}(x + iy)\phi_{n}^{t}(x + iy)) \\
&&\\
& = & - \phi_{n}(x + iy)D_{0}^{-1}\psi^{t}_{n}(x + iy) + \phi_{n}(x + iy)M^{\phi}(x + iy)\phi_{n}^{t}(x + iy).
\end{array}
$$

Note that, for each $n\in\N$, both $\psi_{n}(z)$ and $\phi_{n}(z)$ depend continuously on $z\in\C$, and so, for each $x\in\R$, 
$$
\left\{
\begin{array}{lll}
\lim_{y \downarrow 0} \psi_{n}(x + iy) & = & \psi_{n}(x), \\
& & \\
\lim_{y \downarrow 0} \phi_{n}(x + iy) & = & \phi_{n}(x),
\end{array}\right.
$$
with $\psi_{n}(x), \phi_{n}(x) \in M(l, \mathbb{R})$. Thus, it follows that
$$
\lim_{y \downarrow 0} \Im\left[\phi_{n}(x + iy)D_{0}^{-1}\psi^{t}_{n}(x + iy)\right]  =  0 
$$
and that
$$
\lim_{y \downarrow 0} \Im\left[\phi_{n}(x + iy)M^{\phi}(x + iy)\phi_{n}^{t}(x + iy) \right] = \phi_{n}(x)\lim_{y \downarrow 0}\Im\left[M^{\phi}(x + iy)\right]\phi_{n}^{t}(x), 
$$
if $\lim_{y \downarrow 0}\Im\left[M^{\phi}(x + iy)\right]$ exists. Then, given that $M^\phi(z)$ is a matrix-valued Herglotz function, such limit exists for 
$\kappa$-a.e. $x\in\R$, from which follows that for $\kappa$-a.e. $x\in\R$,
\begin{equation*}\label{MeG}
\lim_{y \downarrow 0} \Im\left[G(n, n; x + iy)\right] = \phi_{n}(x)\lim_{y \downarrow 0} \Im\left[M^{\phi}(x + iy)\right]\phi_{n}^{t}(x).
\end{equation*}
%
\end{proof}

One needs a final ingredient. Let $\mu_{1}$ and $\mu_{2}$ be finite positive Borel measures on $\R$, and set $\min(\mu_{1}, \mu_{2})$ as the measure given by the law
\begin{equation*}
\min(\mu_{1}, \mu_{2})(\Lambda):= \inf\{\mu_{1}(\Lambda_{1}) + \mu_{2}(\Lambda_{2})\mid \Lambda \subseteq \Lambda_{1} \cup \Lambda_{2}\}.
\end{equation*}

\begin{4}
\label{teo.sup.r.ac}
Let, for each $r\in\{1,\ldots,l\}$, $\mathcal{S}_{r}$ be given by \eqref{def.val.r.matriz}. Then, for $\eta_{r} d\kappa$-a.e. $x\in\R$, one has $x \in \mathcal{S}_{r}$, where $\eta_{r} d\kappa = \min (\eta^{\phi}_{r} d\kappa, \eta^{\psi}_{r} d\kappa)$, with  
$$
\eta^{\phi, \psi}_{r}: \dom(\eta^{\phi, \psi}_{r})  \rightarrow  (0, \infty) 
$$
given by
$$
\eta^{\phi, \psi}_{r}(x)  :=  s_{r}\left[\lim_{y \downarrow 0} \Im\left[M^{\phi, \psi}_+(x + iy)\right]\right], 
$$
where
$$
\dom(\eta_r^{\phi, \psi}) :=  \left\{x \in \mathbb{R}\mid 0 < s_{r}\left[\lim_{y \downarrow 0} \Im\left[M^{\phi, \psi}_+(x + iy)\right]\right] < \infty\right\}.
$$ 
\end{4}
\begin{proof}
Let $r\in\{1,\ldots,l\}$. 
It follows from the commutative property of the trace and from Proposition~\ref{Th6.7Haeal}-c (note that $\Im\left[M^{\phi}(x + iy)\right]\ge 0$) that for each $x,y\in\R$ and each $n\ge 0$,
$$
\begin{array}{lll}
 \tr [\phi_{n}(x)\Im\left[M^{\phi}(x + iy)\right]\phi_{n}^{t}(x)] & = &  \tr [\Im\left[M^{\phi}(x + iy)\right]\phi_{n}^{t}(x)\phi_{n}(x)] \\
& & \\
 &=& \sum_{k=1}^ls_k[\Im\left[M^{\phi}(x + iy)\right]\phi_{n}^{t}(x)\phi_{n}(x)]\\
&&\\
 & \geq & \left(s_{r}[\Im\left[M^{\phi}(x + iy)\right]]\right) \left(s_{l - r + 1}[\phi^{t}_{n}(x)\phi_{n}(x)] \right)\\
& & \\
& = & \left(s_{r}[\Im\left[M^{\phi}(x + iy)\right]]\right) \left(s^{2}_{l - r + 1}[\phi_{n}(x)]\right),
\end{array}
$$
and so, by the continuity of the trace and Lemma~\ref{lema.g.ident}, that
\begin{equation}\label{MLEG}
  \lim_{y\downarrow 0}\left(s_{r}[\Im\left[M^{\phi}(x + iy)\right]]\right) \left(s^{2}_{l - r + 1}[\phi_{n}(x)]\right)d\kappa(x)\le \tr\left[\lim_{y\downarrow 0}\Im G^\phi(n,n;x+iy)\right]d\kappa(x).
\end{equation}

One gets from Stieltjes Inversion Formula the identity
$$
\lim_{y \downarrow 0} \frac{1}{\pi} \int  \Im\left[G(n, n; x + iy)\right] d\kappa(x) = \mu_{\textbf{e}_{k, n}, \textbf{e}_{j, n}}(\mathbb{R});
$$
since for each $j,k\in\{1,\ldots,l\}$,
$$
\mu_{\textbf{e}_{k, n}, \textbf{e}_{j, n}}(\mathbb{R})  =  \left\langle \textbf{e}_{k, n}, \textbf{e}_{j, n} \right\rangle =\delta_{k, j},
$$
one concludes that
$$
\frac{1}{\pi} \lim_{y \downarrow 0} \int\Im \left[ G^\phi(n, n; x + iy)_{k, j} \right] d\kappa(x) = \delta_{k, j},
$$
that is, 
\begin{equation}\label{TR1}
\tr\left[\lim_{y \downarrow 0} \int\Im \left[ G^\phi(n, n; x + iy) \right] d\kappa(x)\right] = l\pi.
\end{equation}

By combining relations~\eqref{MLEG} and~\eqref{TR1} and by the definition of the measure $\eta^{\phi}_{r} d\kappa$, one gets 

\begin{equation}
\label{des.sing.phi}
\int s^{2}_{l - r + 1}[\phi_{n}(x)]\eta_r^\phi(x) d\kappa(x) \leq l\pi;
\end{equation}
since the previous argument is also valid for $\psi$ in place of $\phi$, one also gets
\begin{equation}
\label{des.sing.psi}
\int s^{2}_{l - r + 1}[\psi_{n}(x)] \eta^\psi_r(x)d\kappa(x) \leq l\pi.
\end{equation}


Now, since $H^\psi_j$ is a finite rank perturbation of $H^\phi_j$, for each $j\in\{1,\ldots,l\}$ (see the proof of Proposition \ref{prop.supor.pertur}), it follows that the sets $\dom(\eta_r^{\phi})$ and $\dom(\eta_r^{\psi})$ 
are equal up to a set of zero Lebesgue and spectral measures (see the proof of Proposition~3.3 in~\cite{simon99}). Therefore, $\eta_r^{\phi}d\kappa$ and $\eta_r^{\psi}d\kappa$ are equivalent measures, and by the definition of $\eta$ and relations \eqref{des.sing.phi}, \eqref{des.sing.psi}, it follows that for each $n \in \N$, 
\begin{equation}
\label{des.phi.psi.eta}
\int \left( s^{2}_{r - l + 1}[\phi_{n}(x)] + s^{2}_{l - r + 1}[\psi_{n}(x)] \right) \eta_r(x) d\kappa(x)\leq 2l\pi,
\end{equation}
and consequently, that 
$$
\int \frac{1}{L} \sum^{L}_{n = 1} \left( s^{2}_{r - l + 1}[\phi_{n}(x)] + s^{2}_{l - r + 1}[\psi_{n}(x)] \right) \eta_r(x) d\kappa(x)\leq 2l\pi.
$$ 

Finally, one has from Fatou Lemma that 
$$
\int \liminf_{L \rightarrow \infty} \frac{1}{L} \sum^{L }_{n = 1} \left( s^{2}_{l - r + 1}[\phi_{n}(x)] + s^{2}_{r - l + 1}[\psi_{n}(x)] \right) \eta_r(x) d\kappa(x) \leq 2l\pi,
$$
and so, for $\eta_r d\kappa$-a.e. $x\in\R$, that 
$$
\liminf_{L \rightarrow \infty} \frac{1}{L} \sum^{L}_{n =1} \left( s^{2}_{l - r + 1}[\phi_{n}(x)] + s^{2}_{r - l + 1}[\psi_{n}(x)] \right) < \infty.
$$
\end{proof}


\begin{7}
\label{coro.ces.mul.max}
Let $r\in\{1,\ldots,l\}$.  The restriction of the absolutely continuous spectrum of the operator $H^{\phi}_+$ to the set $\mathcal{S}_{r} \setminus \mathcal{S}_{r + 1}$ has multiplicity at most $r$.
\end{7}
\begin{proof}
We begin noting that if $x\in\R$ is such that 
$$
0 < s_{r}[\lim_{y \downarrow 0} \Im[M^{\phi}(x + iy)]] < \infty,
$$
then $x \in \bigcup^{l}_{k = r} \Sigma^{\phi}_{ac, k}$ (with $\Sigma^{\phi}_{ac, k}$ given in the statement of Proposition \ref{porp.sup.ac}).

Thus, it follows from Proposition \ref{teo.sup.r.ac} that if $x \notin\mathcal{S}_{r + 1}$, then
$$
x \notin \bigcup^{l}_{j = r + 1} \Sigma_{ac, j}^{\phi}. 
$$ 
The result follows now from Proposition~\ref{porp.sup.ac}.
\end{proof}

\

\begin{proof2}
  It follows from Theorems \ref{prop.suporte.cesaro} and \ref{coro.ces.mul.max} that $\sigma_{ac,r}(H^\phi_+)=\overline{\mathcal{S}_{r} \setminus \mathcal{S}_{r + 1}}^{ess}$. 
  Now, given that any self-adjoint extension of $H^{+}_{min}$ is a finite rank perturbation of $H^{\phi}_+$, the result is a consequence of the arguments presented in the proof of  Proposition \ref{prop.supor.pertur}. 
\end{proof2}


\section{Constancy of the Absolutely Continuous Spectral Components}
\label{constancy}
\zerarcounters

In the context of dynamically defined scalar Schr\"odinger operators, 
one may highlight the so-called almost-periodic operators, a class of operators that includes the 
 limit-periodic and quasi-periodic operators (see~\cite{damanik17,marx15} for an account of the main results in the theory). In what follows, we extend these definitions to the setting of matrix-valued Jacobi operators.

\begin{1}[Almost-periodic Potential]
\label{def.almost}
A bounded function $V: \mathbb{Z} \rightarrow M(l, \mathbb{C})$ is called an almost-periodic potential if the sequence $(S^{n}(V))_{n \in \mathbb{Z}}$ is a compact set in $l^{\infty}(\mathbb{Z}; M(l, \mathbb{C}))$, where
$$
\begin{array}{lccc}
S: &  l^{\infty}(\mathbb{Z}; M(l, \mathbb{C})) & \rightarrow & l^{\infty}(\mathbb{Z}; M(l, \mathbb{C})) \\
& (V_{n}) & \hookrightarrow  & (V_{n + 1})
\end{array}
$$ 
is the full-shift map.
\end{1}

Let $W: \mathbb{N} \rightarrow M(l, \mathbb{C})$; the hull of $W$ is the set  
\begin{equation}
\label{eq.envoltoria}
\hull(W) := \overline{\{S^{n}(W), n \in \mathbb{N}\}},
\end{equation}
with the closure taken with respect to the topology of uniform convergence in the space $l^{\infty}(\mathbb{Z}; M(l, \mathbb{C}))$. If $W$ is almost-periodic, then $\hull(W)$ is compact. Clearly, every periodic potential is almost-periodic. 

\begin{1}[Limit-periodic Potential]
A potential $V: \mathbb{Z} \rightarrow M(l, \mathbb{C})$ is called limit-periodic if there exists a family of periodic bilateral sequences $W^{(m)}: \mathbb{Z} \rightarrow M(l, \mathbb{C})$ such that
$$
\lim_{m \rightarrow \infty} \sup_{k \in \mathbb{Z}} \left\| V_{k} - W^{(m)}_{k} \right\| = 0.
$$
\end{1}

Note that every limit-periodic potential is also an almost-periodic potential. 
 


Now, given an almost-periodic potential $V$, denote its hull by $\Omega$. Then, $\Omega$ is a compact topological subspace of  $l^{\infty}(\mathbb{Z}; M(l, \mathbb{C}))$ and one can endow $\Omega$ with an abelian group structure. Namely, if $W_{1}, W_{2}$ are potentials such that there exist $k_{1}, k_{2} \in \mathbb{Z}$ so that $W_{1} = S^{k_{1}}(V)$ and $W_{2} = S^{k_{2}}(V)$, then one may define the product of $W_1$ and $W_2$ by the law
$$
W_{1} \ast W_{2} := S^{k_{1} + k_{1}}(V).
$$ 

One may also extend such product to $\Omega$, which becomes  a compact topological group, and so, endowed with a Haar measure.

Therefore, for each almost-periodic potential $V$, there exist a continuous function $f: \Omega \rightarrow  M(l, \mathbb{R})$, a minimal translation $T: \Omega \rightarrow \Omega$ and $\omega\in\Omega$ such that, for each $n\in\Z$,
\begin{equation}
\label{eq.almost.prot}
V_{n} = V_{\omega}(n) = f(T^{n}\omega).
\end{equation}

One may extend the previous result and define the notion of a minimal matrix-valued Jacobi operator. Let $\Omega \neq \emptyset$ be a compact metric space, let $T: \Omega \rightarrow \Omega$ be a minimal invertible map and let 
$D, V: \Omega \rightarrow S(l, \mathbb{R})$, $l\in\N$, be continuous maps such that, for each $\omega\in\Omega$, $\det (D(\omega))\neq 0$ and $D(\omega)$, $V(\omega)$ are both symmetric. A minimal operator is, then, a family of operators in $(\mathbb{C}^{l})^{\mathbb{Z}}$ given, for each $\omega \in \Omega$, by \eqref{eq.ope.din}, that is, for each $n\in\Z$ one has
\[[H_{\omega} \textbf{u}]_{n} := D(T^{n - 1}\omega) \textbf{u}_{n - 1} + D(T^{n}\omega) \textbf{u}_{n + 1} + V(T^{n}\omega) \textbf{u}_{n}.
\]

One of the consequences of Theorems \ref{prop.suporte.cesaro} and \ref{teo.sup.r.ac}, as discussed  in \cite{simon99} for the scalar  case, is the constancy of the absolutely continuous spectral components for minimal operators. 

The concept of minimal operator has some topological consequences.
 
\begin{1}[Right Limit and Left Limit of a Potential]
Let $W: \mathbb{Z} \rightarrow M(l, \mathbb{R})$ be an almost-periodic potential. One says that $W_{\omega_2}$ is the right limit of $W_{\omega_1}$ if there exists a sequence $(n_{j})_{j\ge 1}$ with $\lim_{j \rightarrow \infty}n_{j} = \infty$ such that, for each $n \in \mathbb{Z}$, one has
$$
\lim_{j \rightarrow \infty} W_{\omega_1}(n - n_{j}) = W_{\omega_2}(n).
$$
In the same conditions, one says that $W_{\omega_2}$ is the left limit of $W_{\omega_1}$ if, for each $n \in \mathbb{Z}$,
$$
\lim_{j \rightarrow \infty} W_{\omega_1}(n + n_{j}) = W_{\omega_2}(n).
$$
\end{1}

Note that if, for a fixed $\omega\in\Omega$, $V_{\omega}$ is the potential associated with a minimal operator and $W \in \hull(V_{\omega})$, then, by definition, $W$ is the left limit of $V_\omega$. 


\begin{4}
\label{prop.const.almost}
Let $\omega_1,\omega_2\in\Omega$ and let $H_{\omega_1}$, $H_{\omega_2}$ be the matrix-valued dynamically defined Jacobi operators given by \eqref{eq.ope.din}, respectively. 
If $D_{\omega_2}$ is the left limit of $D_{\omega_1}$ and if $V_{\omega_2}$ is the left limit of $V_{\omega_1}$, then, for each $k \in \{1,2,\ldots,l\}$, 
$$
\sigma_{ac, k}(H_{\omega_1}^{\phi}) = \sigma_{ac,k}(H_{\omega_2}^{\phi}),
$$ 
where $\sigma_{ac,k}(H_{\omega_j}^{\phi})$ is the absolutely continuous spectrum of multiplicity $k$ of the operator $H_{\omega_j,+}^{\phi}$.
\end{4}
\begin{proof}
  Let $k\in\{1,\ldots,l\}$. 
  By hypothesis, there exists a sequence $(n_{j})_{j\ge 1}$ with $\lim_{j \rightarrow \infty}n_{j} = \infty$ such that, for each $m \in \N$, one has
$$
\begin{array}{lll}
\lim_{j \rightarrow \infty} D_{\omega_1}(m + n_{j}) & = & D_{\omega_2}(m), \\
\lim_{j \rightarrow \infty} V_{\omega_1}(m + n_{j}) & = & V_{\omega_2}(m).
\end{array}
$$

Let $\phi^{\omega_1}$ and $\psi^{\omega_1}$ be, respectively, the Dirichlet and Neumann solutions to the eigenvalue equation~\eqref{eq.aut.Jacobi} for the operator $H_{\omega_1}$, and denote $D_\omega(m)$ (respectively, $V_\omega(m)$) by $D^\omega_m$ (respectively, $V^\omega_m$). Then, for each 
$j,m\in \N$ and each $x\in\R$, one has
\begin{eqnarray}
\label{eq.lim.esq.phi}
\left\{
\begin{array}{lll}
  D^{\omega_1}_{m + n_{j}}\phi^{\omega_1}_{m + n_{j} + 1} + D^{\omega_1}_{m + n_{j} - 1}\phi^{\omega_1}_{m + n_{j} - 1}
+ V^{\omega_1}_{m + n_{j}}\phi^{\omega_1}_{m + n_{j}} & = & x \phi^{\omega_1}_{m + n_{j}}, \\
& & \\
D^{\omega_1}_{m + n_{j}}\psi^{\omega_1}_{m + n_{j} + 1} + D^{\omega_1}_{m + n_{j} - 1}\psi^{\omega_1}_{m + n_{j} - 1}
+ V^{\omega_1}_{m + n_{j}}\psi^{\omega_1}_{m + n_{j}} & = & x \psi^{\omega_1}_{m + n_{j}}.
\end{array}\right.
\end{eqnarray}

It follows from the compacity of $\Omega$ and from the continuity of the solutions to the eigenvalue equation~\eqref{eq.aut.Jacobi}  with respect to $\omega\in\Omega$ that (by taking a subsequence if necessary) for each $m \in \N$, the limits  
$$
\begin{array}{lll}
B_{m} & := & \lim_{j  \rightarrow \infty} \phi^{\omega_1}_{m + n_{j}}, \\
& & \\
C_{m} & := & \lim_{j  \rightarrow \infty} \psi^{\omega_1}_{m + n_{j}}
\end{array}
$$
do exist. Now, by letting $j\rightarrow\infty$ in both members of the identities \eqref{eq.lim.esq.phi}, one gets for each $m\in\N$ that
$$
\left\{
\begin{array}{lll}
D^{\omega_2}_{m}B_{m + 1} + D^{\omega_2}_{m - 1}B_{m - 1} + V^{\omega_2}_{m}B_{m} & = & x B_{m}, \\
& & \\
D^{\omega_2}_{m}C_{m + 1} + D^{\omega_2}_{m - 1}C_{m - 1} + V^{\omega_2}_{m}C_{m} & = & x C_{m}.
\end{array}\right.
$$
Therefore, one concludes that $(B_{n})_{n \in \N}$ and $(C_{n})_{n \in \N}$ are solutions to the eigenvalue equation~\eqref{eq.aut.Jacobi} at $x$ for the operator $H_{\omega_2}$. 

Now, as in the proof of Lemma~2.4 in~\cite{vieiracarvalho21}, by applying relation \eqref{wronski.matriz} to $\phi^{\omega_1}$ and $\psi^{\omega_1}$, one gets for each $j,m \in \N$,
$$
\left\{
\begin{array}{lll}
(\psi^{\omega_1}_{m + n_{j}})^{t}D^{\omega_1}_{m + n_{j}}\psi^{\omega_1}_{m + n_{j} + 1} - (\psi^{\omega_1}_{m + n_{j} + 1})^{t}D^{\omega_1}_{m + n_{j}}\psi^{\omega_1}_{m + n_{j}} & = & 0, \\ 
&&\\
(\psi_{m + n_{j}}^{\omega_1})^{t}D^{\omega_1}_{m + n_{j}}\phi^{\omega_1}_{m + n_{j} + 1} - (\psi^{\omega_1}_{m + n_{j} + 1})^{t}D^{\omega_1}_{m + n_{j}}\phi^{\omega_1}_{m + n_{j}} & = & D^{\omega_1}_{0}, \\
&&\\
(\phi^{\omega_1}_{m + n_{j}})^{t}D^{\omega_1}_{m + n_{j}}\psi^{\omega_1}_{m + n_{j} + 1} - (\phi^{\omega_1}_{m + n_{j} + 1})^{t}D^{\omega_1}_{m + n_{j}}\psi^{\omega_1}_{m + n_{j}} & = & - D^{\omega_1}_{0}, \\
&&\\
(\phi^{\omega_1}_{m + n_{j}})^{t}D^{\omega_1}_{m + n_{j}}\phi^{\omega_1}_{m + n_{j} + 1} - (\phi^{\omega_1}_{m + n_{j} + 1})^{t}D^{\omega_1}_{m + n_{j}}\phi^{\omega_1}_{m + n_{j}} & = & 0,
\end{array}\right.
$$
and then, by letting $j\rightarrow\infty$ in both members of these identities, it follows that for each $m\in\N$,
$$
\left\{
\begin{array}{lll}
(C_{m})^{t}D^{\omega_2}_{m}C_{m + 1} - (C_{m + 1})^{t}D^{\omega_2}_{m}C_{m} & = & 0, \\ 
&&\\
(C_{m})^{t}D^{\omega_2}_{m}B_{m + 1} - (C_{m + 1})^{t}D^{\omega_2}_{m}B_{m} & = & D^{\omega_1}_{0}, \\
&&\\
(B_{m})^{t}D^{\omega_2}_{m}C_{m + 1} - (B_{m + 1})^{t}D^{\omega_2}_{m}C_{m} & = & - D^{\omega_1}_{0}, \\
&&\\
(B_{m})^{t}D^{\omega_2}_{m}B_{m + 1} - (B_{m + 1})^{t}D^{\omega_2}_{m}B_{m} & = & 0.
\end{array}\right.
$$

Hence, the sequences $(B_{n})_{n \in \N}$ and $(C_{n})_{n \in \N}$ satisfy the same relations used in Section \ref{sec.mat.tra} for obtaining the essential supports \eqref{def.val.r.matriz} with respect to $\psi$ and $\phi$. Therefore, one concludes that $B$ and $C$ constitute a basis for the space of the solutions to the eigenvalue equation~\eqref{eq.aut.Jacobi} for $H_{\omega_2}$. Then, by Theorem~\ref{coro.sr}, if one defines for $r\in\{1,\ldots,l\}$ the sets 
$$
\mathcal{T}_{r}^{\omega_2} := \{x \in \mathbb{R}\mid \liminf_{L \rightarrow \infty} \frac{1}{L} \sum^{L}_{n = 1} s^{2}_{l - r + 1}[B_{n}(x)] + s^{2}_{l - r + 1}[C_{n}(x)] < \infty \},
$$
then $\overline{\mathcal{T}_{r}^{\omega_2} \setminus \mathcal{T}_{r + 1}^{\omega_2}}^{ess}$ corresponds to the absolutely continuous spectral component of multiplicity $r$ of any self-adjoint extension of the operator $H_{\omega_2,+}^{min}$, and in particular, of $H_{\omega_2,+}^{\phi}$.

Now, as in the proof of Theorem \ref{teo.sup.r.ac} (inequality \eqref{des.phi.psi.eta}), one has for each $n \in \N$, 
$$
\dint \left( s^{2}_{l - k + 1}[\phi^{\omega_1}_{n}(x)] + s^{2}_{l - k + 1}[\psi^{\omega_1}_{n}(x)] \right) \eta_k^{\omega_1}(x) d\kappa(x) \leq 2l\pi,
$$
with $\eta_k^{\omega_1} d\kappa$ supported in $\sigma_{ac, k}(H^{\phi}_{\omega_1})$; consequently, it follows that for each $L, j \in \N$,
$$
\begin{array}{llll}
& \dint \dfrac{1}{L} \sum^{L}_{n = 1}
\left( s^{2}_{l - k + 1}[\phi^{\omega_1}_{n + n_{j}}(x)] + s^{2}_{l - k + 1}[\psi^{\omega_1}_{n + n_{j}}(x)] \right) \eta_k^{\omega_1}(x) d\kappa(x) & = & \\
& & & \\
= & \dint \dfrac{1}{L} \sum^{n_{j} + L}_{m = n_{j} + 1}
\left( s^{2}_{l - k + 1}[\phi^{\omega_1}_{m}(x)] + s^{2}_{l - k + 1}[\psi^{\omega_1}_{m}(x)] \right) \eta_k^{\omega_1}(x) d\kappa(x) & \leq & 2l\pi.
\end{array}
$$
Then, it follows from Dominated Convergence Theorem that 
$$
\int \frac{1}{L} \sum^{L}_{n = 1}
\left( s^{2}_{l - k + 1}[B_{n}(x)] + s^{2}_{l - k + 1}[C_{n}(x)] \right) \eta_k^{\omega_1}(x) d\kappa(x) \leq 2l\pi,
$$
and, finally, from Fatou Lemma that
\begin{equation}
\label{eq.almost.per.cesaro}
\liminf_{L \rightarrow \infty} \frac{1}{L} \sum^{L}_{n = 1} s^{2}_{l - k + 1}[B_{n}(x)] + s^{2}_{l - k + 1}[C_{n}(x)]  < \infty,
\end{equation}
for a.e. $x\in\R$ with respect to the measure $\eta^{\omega_1}_kd\kappa$.

Hence, by Theorem \ref{prop.suporte.cesaro} and by the remark about $\mathcal{T}_k$ presented above, the set inclusion $\mathcal{T}_k^{\omega_1}\subset\mathcal{T}_k^{\omega_2}$ follows up to sets of zero Lebesgue and $\eta_k^{\omega_1}d\kappa$ measures (namely, it follows from Theorem~\ref{prop.suporte.cesaro}, Proposition~\ref{teo.sup.r.ac} and the aforementioned remark that $\mathcal{T}_k^{\omega_1}$ is a minimal support for $\eta^{\omega_1}_kd\kappa$). 
Naturally, the set inclusion $\mathcal{T}_k^{\omega_2}\subset\mathcal{T}_k^{\omega_1}$, valid up to sets of zero Lebesgue and $\eta_k^{\omega_2}d\kappa$ measures, 
%
can be obtained using the same arguments as above. Therefore,
\[\sigma_{ac, k}(H_{\omega_1}^{\phi})=\overline{\mathcal{T}_k^{\omega_1}\setminus\mathcal{T}_{k+1}^{\omega_1}}^{ess}=\overline{\mathcal{T}_k^{\omega_2}\setminus\mathcal{T}_{k+1}^{\omega_2}}^{ess}=\sigma_{ac, k}(H_{\omega_2}^{\phi}),\]
and we are done.
\end{proof}


\begin{7}
\label{teo.cont.ac.almost}
Let $(H_{\omega})_{\omega \in \Omega}$ be the family of minimal matrix-Jacobi operators given by~\eqref{eq.ope.din}. Then, for each $\omega_{0}, \omega_{1} \in \Omega$ and each $k\in\{1,\ldots,l\}$, one has 
$$
\sigma_{ac, k}(H_{\omega_{1}}^{\phi}) = \sigma_{ac, k}(H_{\omega_{0}}^{\phi}),
$$
where $H_{\omega}^{\phi}$ and $\sigma_{ac,k}(H_{\omega}^{\phi})$ are defined as in Proposition \ref{prop.const.almost}.
\end{7}

\begin{proof}
Let, for each $\omega \in \Omega$, $D^{\omega}$ denote the bilateral sequence $(D_{\omega}(n))_{n \in \mathbb{Z}}$. It follows from the minimality of $T$ that for each $\omega_{0}, \omega_{1} \in \Omega$, $D^{\omega_{1}} \in \hull(D^{\omega_{0}})$. Hence, $D^{\omega_{0}}$ is the right limit of $D^{\omega_{1}}$. In the same way, $V^{\omega_{0}}$ is the right limit of $V^{\omega_{1}}$. The result follows now from Proposition~\ref{prop.const.almost}.
\end{proof}


\

\begin{proof3}
  Since $\omega_{1} \in \overline{\mathcal{O}(\omega_{0})}$, where $\mathcal{O}(\omega_0):=\{T^n\omega_0\mid n\in\N\}$, there exists a sequence $(n_{j})_{j\in\N}$ with $n_j\rightarrow\infty$ such that $\lim_jT^{n_{j}}\omega_{1} = \omega_{0}$. In this way, if one denotes by $\psi^{\omega_{0}}$ and $\phi^{\omega_{0}}$ the Neumann and Dirichlet solutions to the eigenvalue equation for the operator $H_{\omega_{0}}$, one can prove (as in the proof of Proposition \ref{prop.const.almost}) that there exist sequences $(B_{n})_{n \in \mathbb{Z}_-}$ and $(C_{n})_{n \in \mathbb{Z}_-}$ of solutions to the eigenvalue equation~\eqref{eq.aut.Jacobi} for $H_{\omega_{1}}$ such that for each $m \in \mathbb{Z}_-$,
$$
\begin{array}{lll}
B_{m} & := & \lim_{j  \rightarrow \infty} \phi^{\omega_{0}}_{m + n_{j}}, \\
& & \\
C_{m} & := & \lim_{j  \rightarrow \infty} \psi^{\omega_{0}}_{m + n_{j}}.
\end{array}
$$

In particular, for each $r\in\{1,\ldots,l\}$ and each $L \in  \mathbb{N}$, one has
$$
\frac{1}{L} \sum^{-L}_{m = -1} s^{2}_{l - r + 1}[B_{m}(x)] + s^{2}_{l - r + 1}[C_{m}(x)] 
=
\lim_{j  \rightarrow \infty} 
\frac{1}{L} \sum^{-L}_{m = -1} s^{2}_{l - r + 1}[\phi^{\omega_{0}}_{m + n_{j}}(x)] + s^{2}_{l - r + 1}[\psi^{\omega_{0}}_{m + n_{j}}(x)]. 
$$

Now, it follows from the same arguments presented in the proof of Proposition \ref{prop.const.almost} that for each $m \in \mathbb{Z}_-$ and each $n_{j}$, 
$$
\int \left( \frac{1}{L} \sum^{-L}_{m = -1} s^{2}_{l - k + 1}[\phi^{\omega_{0}}_{m + n_{j}}(x)] + s^{2}_{l - k + 1}[\psi^{\omega_{0}}_{m + n_{j}}(x)] \right) \eta^{\omega_{0}}_k(x) d\kappa(x) \leq 2 \pi l,
$$
where $\eta_k^{\omega_{0}} d\kappa$ is supported in $\sigma_{ac, k}((H^{\phi}_{+})_{\omega_{0}})$. By Dominated Convergence Theorem, one gets 
$$
\int \left( \frac{1}{L} \sum^{-L}_{m = -1} s^{2}_{l - k + 1}[B_{m}(x)] + s^{2}_{l - k + 1}[C_{m + n_{j}}(x)] \right) \eta_k^{\omega_{0}}(x) d\kappa(x) \leq 2 \pi l, 
$$
and then, by Fatou Lemma, it follows that
$$
\liminf_{L \rightarrow \infty} \frac{1}{L} \sum^{-L}_{m = -1} s^{2}_{l - k + 1}[B_{m}(x)] + s^{2}_{l - k + 1}[C_{m + n_{j}}(x)] < \infty 
$$
is true for a.e $x\in\R$ with respect to the measure $\eta_k^{\omega_{0}} d\kappa$. Therefore, it follows from 
the arguments presented in the proof of Proposition  \ref{prop.const.almost} that 
$$
\sigma_{ac, k}((H^{\phi}_{-})_{\omega_{0}})=\sigma_{ac, k}((H^{\phi}_{+})_{\omega_{1}}),
$$ 
and also by 
Proposition \ref{prop.const.almost} that 
$$
\sigma_{ac, k}((H^{\phi}_{-})_{\omega_{0}}) = \sigma_{ac, k}((H^{\phi}_{-})_{\omega_{1}})
$$ 
and
$$
\sigma_{ac, k}((H^{\phi}_{+})_{\omega_{0}}) = \sigma_{ac, k}((H^{\phi}_{+})_{\omega_{1}}).
$$

The result is now a consequence of Proposition \ref{prop.supor.pertur.z}. Namely, it follows from the previous relations that for each $k\in\{1,\ldots,l\}$ and each $\omega_0,\omega_1\in\Omega$,
\[\sigma_{ac, k}((H^{\phi}_{-})_{\omega_0})=\sigma_{ac, k}((H^{\phi}_{+})_{\omega_0})=\sigma_{ac, k}((H^{\phi}_{-})_{\omega_1})=\sigma_{ac, k}((H^{\phi}_{+})_{\omega_1}),\]
and so, one has from Proposition~\ref{prop.supor.pertur.z} that
\[\sigma_{ac, 2k}(H_{\omega_0})=\sigma_{ac, 2k}((H^{\phi}_{+})_{\omega_0}\oplus(H^{\phi}_{-})_{\omega_0})=\sigma_{ac, 2k}((H^{\phi}_{+})_{\omega_1}\oplus(H^{\phi}_{-})_{\omega_1})=\sigma_{ac, 2k}(H_{\omega_1})\]
and that
\[\sigma_{ac, 2k+1}(H_{\omega_0})=\emptyset=\sigma_{ac, 2k+1}(H_{\omega_1}).\]
\end{proof3}

\section{Appendix}

Here, we present the results involving the Frobenius norm and the singular values of square matrices with complex entries that are used in this work.

A matrix $A\in M (l,\C)$ is said to be Hermitean if $A^\ast=A$. A Hermitean matrix $A$ is said to be semidefinite positive if for each $\mathbf{u}\in \C^l\setminus\{\mathbf{0}\}$, $\mathbf{u}A\mathbf{u}^t\ge 0$. One represents this fact as $A\ge0$. Note that for each $A\in M(l,\C)$, $A^\ast A$ is hermitean and positive semidefinite.

One may also define a partial order in the set of positive semidefinite matrices. Namely, if $A,B\in M(l,\C)$ are positive semidefinite matrices, one says that $A\ge B$ if $\mathbf{u}A\mathbf{u}^t\ge\mathbf{u} B\mathbf{u}^t$. In other words, $A\ge B$ if, and only if, $A-B\ge 0$.

If $A$ is a positive semidefinite matrix, there exists only one positive semidefinite matrix $B$ such that $B^2=A$, the so-called square root of $A$: $B=\sqrt{A}$.

One may easily check that the map $\langle\cdot,\cdot\rangle:M(l,\C)\times M(l,\C)\rightarrow\C$ given by the law
\[\langle A,B\rangle=\tr[A^\ast B]\]
is an inner product in $M(l,\C)$.

\begin{1}[Frobenius norm]
\label{frobeniusnorm}  
  The norm induced in $M(l,\C)$ by the previous inner product is called Frobenius norm:
  \[\Vert A\Vert_F:=\sqrt{\tr[A^\ast A]}.\]
\end{1}
  
\begin{1}[Singular value of a square matrix]\label{sedvalsing}
  The singular values of $A\in M(l,\C)$ are defined as the eigenvalues of $\sqrt{A^\ast A}$. We denote them by
\[s_1 [A]\ge s_2 [A]\ge\ldots\ge s_{l-1} [A]\ge s_l [A],\]
with $s_k[A] =\lambda_k(\sqrt{A^\ast A})\ge 0$, including the multiplicity.
\end{1}

\begin{4}[Singular values and norms; Theorem~6.7 in~\cite{petz2014}]
\label{traceFrobenius}  
Let $A\in M (l,\C)$. Then,
\[\Vert A\Vert = s_1[A],\]
where $\Vert A\Vert:=\sup_{\Vert\mathbf{u}\Vert=1}\Vert A\mathbf{u}\Vert$ stands for the operator norm of $A$.

On the other hand,
\[\Vert A\Vert_F=\sqrt{s_1^2[A] + s_2^2[A] +\ldots+ s_l^2[A]}.\]
\end{4}

\begin{4}[Theorem~6.7 in~\cite{petz2014}]
\label{Th6.7Haeal}
Let $A, B\in M (l,\C)$.
\begin{itemize}
\item[(a)]
If $0\le A\le B$, then, for each $k\in\{1,\ldots,l\}$,
\[s_k [A]\le s_k [B].\]
\item[(b)] If $k + m - 1\le l$, then 
\[s_{k+m-1}[A + B]\le s_k[A] + s_m[B].\]
\item[(c)] For each $k\in\{1,\ldots,l\}$, one has
  \[\sum_{i=1}^ks_i[AB]\ge\sum_{i=1}^ks_i[A]s_{l-i+1}[B].\]
\item[(d)] If $k,m\in\{1,\ldots,l\}$ are such that $k + m-1\le l$, then
  \[s_l[A]s_m[B]\le s_{k+m-1}[AB]\le s_1[A]s_m[B].\]
\item[(e)] (Minimax expression) If $A\ge 0$, then for each $k\in\{1,\ldots,l\}$,
  \[s_k[A] = \min\{\max\{\langle\mathbf{u},A\mathbf{u}\rangle\mid\mathbf{u}\in\mathcal{M}^\perp,\Vert\mathbf{u}\Vert=1\}\mid\mathcal{M}\;\;\textrm{is a subspace of dimension}\;\; k\}.\]  
\end{itemize}
\end{4}

\begin{remark}
  Proposition~\ref{Th6.7Haeal}-c is actually proven in~\cite{wang97}.
\end{remark}

\begin{center} \Large{Acknowledgments} 
\end{center}
\addcontentsline{toc}{section}{Acknowledgments}

Silas L. Carvalho thanks the partial support by FAPEMIG (Minas Gerais state agency; Universal Project under contract 001/17/CEX-APQ-00352-17) and Fabr\'icio Oliveira thanks the partial support by CAPES (a Brazilian government agency).

We also thanks the anonymous referee, Marcin Moszy\'nski and Grzegorz \'Swiderski for his invaluable suggestions that helped to improve the quality of the manuscript. 

\bibliography{bibfile}{}
\bibliographystyle{acm}

\end{document}